\documentclass[4paper]{article}

\usepackage[nohead,margin=1.0in]{geometry}
\usepackage{graphicx}
\usepackage{amsmath}
\usepackage{amsthm}
\usepackage{amstext}
\usepackage{amsfonts}
\usepackage{amssymb}
\usepackage{amsxtra}
\usepackage{bbm}
\usepackage{url}
\usepackage{mathtools}
\usepackage{booktabs}       
\usepackage{tikz}

\usetikzlibrary{intersections} 
\tikzset{>=latex}
\usepackage{tikzscale}
\usepackage[ruled, vlined,linesnumbered]{algorithm2e}
\usepackage{setspace}
\usepackage{lipsum}		
\usepackage{footnote}
\makesavenoteenv{tabular}
\makesavenoteenv{table}
\usepackage{pgfplots}
\usepgfplotslibrary{fillbetween}

\usepackage{siunitx}
\newtheorem{theorem}{Theorem}[section]
\newtheorem{lemma}[theorem]{Lemma}

\def\RR{\mathbb{R}}

\def\CC{\mathbb{C}}

\newcommand{\norm}[1]{\left\lVert#1\right\rVert}

\newcommand{\abs}[1]{\left| #1 \right|}

\newcommand{\pen}{f}

\newcommand{\prox}{\textrm{prox}}
\DeclareMathOperator*{\argmin}{arg\,min}

\newcommand{\nng}{{\mbox{\scriptsize NNG}}}

\makeatletter
\newcommand\notsotiny{\@setfontsize\notsotiny\@vipt\@viipt}
\makeatother

\def\epsrel{\epsilon_{r}}

\usepackage{cite}

  \usepackage[caption=false,font=footnotesize]{subfig}

\begin{document}

\title{A Wavelet Based Sparse Row-Action Method for Image Reconstruction in Magnetic Particle Imaging}
\date{}

\author{
  Florian~Lieb\thanks{Department of Computer Science, TH Aschaffenburg, 63741 Aschaffenburg, Germany} $\ $and~Tobias~Knopp\thanks{Section for Biomedical Imaging, University Medical Center Hamburg-Eppendorf and Institute for Biomedical Imaging, Hamburg University of Technology, Germany}
}

\maketitle

\begin{abstract}
Magnetic Particle Imaging (MPI) is a preclinical imaging technique capable of  visualizing the spatio-temporal distribution of magnetic nanoparticles. 
The image reconstruction of this fast and dynamic process relies on efficiently solving an ill-posed inverse problem.
Current approaches to reconstruct the tracer concentration from its measurements are either adapted to image characteristics of MPI but suffer from higher computational complexity and  slower convergence or are fast but lack in the image quality of the reconstructed images.
In this work we propose a novel MPI reconstruction method to combine the advantages of both approaches into a single algorithm. 
The underlying sparsity prior is based on an undecimated wavelet transform and is integrated into a fast row-action framework to solve the corresponding MPI minimization problem.
Its performance is numerically evaluated against a classical FISTA approach on simulated and real MPI data. We also compare the results to the state-of-the-art MPI reconstruction methods. 
In all cases, our approach shows better reconstruction results and at the same time accelerates the convergence rate of the underlying row-action algorithm. 
\end{abstract}

\section{Introduction}
When Magnetic Particle Imaging (MPI) was invented in 2005 by Weizenecker and Gleich, a novel research area in medical imaging has emerged alongside computed tomography (CT) and magnetic resonance imaging (MRI) \cite{Gleich2005}.
MPI images the spatial and temporal distribution of superparamagnetic iron oxide nanoparticles (SPIONs) by exploiting the characteristic non-linear magnetization of the particles when external magnetic fields are applied. 
These fields elicit a spatially localized field-free point (FFP) where the dynamic change of the particles magnetic moment results in an induced voltage. 
This voltage, in turn, can be captured by dedicated receive coils \cite{Panagiotopoulos2015}.
This not only takes the imaging process into the sub-millimeter resolution range, but also allows a fast acquisition process enabling real time or interventional imaging \cite{Bakenecker2018}. 
MPI has a diverse range of applications such as vascular imaging \cite{Vaalma2017}, tumor imaging \cite{Yu2017}, stem cell tracking \cite{Bakenecker2018}, hyperthermia therapy \cite{Tay2018} or drug targeting \cite{Le2017}. 

The correlation between the particle position and the induced voltage is encoded in the system matrix, a collection of calibration scans containing the particle behaviour at all discretized positions in the scanner, the so-called field of view (FOV). 
This correlation depends strongly on the chosen trajectory of the FFP where the Lissajous trajectory is one of the most common ones used for real-time imaging \cite{Knopp2008}. 
The system matrix is usually acquired by measuring the system response of a delta sample shifted to each voxel in the FOV. 
A more sophisticated approach avoids mechanical movement of the delta sample and rather shifts the spatial position by utilizing magnetic offset fields. 
Such hybrid system matrix approaches can cut acquisition times down by a factor of 1000 \cite{Graeser2017}. We note that also the $x$-space reconstruction approaches \cite{ozaslan2019fully} can be formulated in terms of a linear system of equations and thus, the reconstruction framework we propose in this work is completely scanner and imaging sequence agnostic. 

With the system matrix $A$, the particle concentration $x$ and the measured voltage $b$ the reconstruction problem amounts to solving the linear system 
\begin{equation}
Ax = b.
\label{eq:scu}
\end{equation}
The system matrix $A\in\CC^{m\times n}$ is usually considered in frequency domain allowing a frequency dependent preprocessing \cite{Knopp2010}. This preprocessing removes noisy rows of the system matrix with a signal-to-noise-ratio (SNR) lower than some predefined threshold. In addition to the SNR-based thresholding, a band-pass filtering is applied to suppress lower frequencies, which are distorted from the excitation frequency \cite{Kluth2019}.

The reconstruction of the tracer concentration is an ill-posed inverse problem. There are currently two different optimization models to solve \eqref{eq:scu} in the context of MPI. The first one is based on regularizing the total energy of the tracer concentration $x$, resulting in the Tikhonov regularized least-squares problem
\begin{equation}
\min_{x\geq0} \frac{1}{2}\norm{Ax-b}_2^2 + \rho\norm{x}_2^2 \qquad \rho>0.
\label{eq:minprob}
\end{equation}
This problem can be efficiently solved using the regularized Kaczmarz approach proposed in \cite{Knopp2010}. As the rows of the system matrix are nearly orthogonal, it shows an initially fast convergence and is therefore frequently used as the state-of-the-art solver in most published papers with only a few numbers of iterations, see for example \cite{Bakenecker2018,Knopp2016a,Moeddel2019}. On the downside, however, Tikhonov regularization results in smoothed edges, limited noise suppression as well as reduced contrast \cite{Storath2017}. In order to overcome these issues, Storath et al.~proposed the second optimization model, i.e., the non-negative Fused Lasso. Instead of the $\ell_2$-prior in \eqref{eq:minprob}, a combination of $\ell_1$- and total variation (TV) regularization priors were evaluated in \cite{Storath2017,ilbey2017}, resulting in the following optimization problem
\begin{equation}
\min_{x\geq0} \frac{1}{2}\norm{Ax - b}^2_2 + \lambda_1\norm{x}_1 + \lambda_2\norm{x}_{TV} \quad \lambda_1,\lambda_2>0.
\label{eq:fl}
\end{equation}
This is motivated by typical statistics of MPI images: a sparse spatial tracer concentration and a sparse image gradient preserving edges. Solving \eqref{eq:fl} in the context of MPI has been shown to significantly improve the reconstructed tracer distribution, see \cite{Storath2017} and \cite{ilbey2017}. In both works different algorithmic approaches for tackling \eqref{eq:fl} are considered, both leading to similar image quality as well as computation time.
The improved image quality that results from the total variation approaches in \cite{Storath2017, ilbey2017} have one major disadvantage, which has hampered a more widespread use of these reconstruction methods up to now: a higher computational complexity paired with a slower convergence rate compared to the classical Kaczmarz approach. In the numerical results section of \cite{Storath2017} this is quantified to be a factor of seven for a typical 3D real data set, i.e., the regularized Kaczmarz is seven times faster than the Fused Lasso approach. Apart from computational aspects, another disadvantage of the Fused Lasso approach is the fact that there are now two regularization parameters, $\lambda_1$ and $\lambda_2$ which impact the outcome of the reconstruction quality. Fine tuning of both parameters simultaneously is then quite challenging.

In many applications of MPI an offline reconstruction of the tracer distribution is sufficient, and hence longer reconstruction times might be tolerated. Yet, some applications, for example in interventional radiology where an immediate visual response for catheter navigation is crucial, require a framework with a feasible real-time reconstruction algorithm \cite{Knopp2016a}. In such cases the regularized Kaczmarz algorithm is still the preferred choice \cite{Moeddel2019}. Unsurprisingly, the larger the data sets get, in particular in three dimensions, a fast and at the same time accurate reconstruction is fundamental.

\textit{Contributions:} In this work we introduce a simple and fast row-action framework in order to reconstruct the tracer concentration in Magnetic Particle Imaging. The novelty lies in combining the basic Kaczmarz based approach from \eqref{eq:minprob} with sparsity inducing priors based on sparse representations of the particle concentration in a discrete wavelet basis. This links the speed of convergence of the underlying row-action method when solving \eqref{eq:minprob} with a single sparsity prior, which is similar to the Fused Lasso approach in \eqref{eq:fl}. Our numerical experiments confirm that the proposed framework leads to a fast and accurate algorithm that performs well and is significantly faster than current state-of-the-art MPI optimization algorithms. The numerical evaluations are based on simulated MPI data with a measured system matrix as well as real three dimensional data provided by the OpenMPIData initiative \cite{Knopp2019}. A Matlab implementation and documentation of the proposed method for reconstructing MPI images is available online \url{https://github.com/flieb/SparseKaczmarzMPI}. 
 
\textit{Outline of the Paper: } In the following section, we first introduce the basic mathematical concepts and notations. In Section \ref{sec:propAlg} we draw the connection between the proposed optimization problem and the MPI problem and introduce our new algorithm. The numerical evaluations are presented in Section \ref{sec:res}, demonstrating the effectiveness of the proposed method on simulated and real MPI data.

\section{Mathematical Background and Notations}

\subsection{Proximal operators}
Proximal operators are introduced as a generalization of convex projections operators and are a popular tool in modern inverse problems allowing the incorporation of sophisticated regularization priors into the optimization problem. 

Let $\pen: \CC^N \rightarrow \RR \cup \{+\infty\}$ be a lower semi-continuous, convex function. For all $x\in\CC^N$, the proximal operator $\prox_{\pen} : \CC^N \rightarrow \CC^N$ $\pen$ is uniquely defined by 
\begin{equation}
\prox_{\pen}(x) = \argmin_{z\in\CC^N} \pen(z) + \frac{1}{2}\norm{x-z}^2_2.
\label{eq:proxop}
\end{equation}

An important result on the composition of $\pen$ with a bounded affine operator $\Phi$ is introduced in \cite[Lemma 2]{Fadili2009} and briefly summarized in the following.
\begin{lemma}\label{lem:proxA}
Let $\Phi:\CC^M \rightarrow \CC^N$ be a linear mapping, such that  $\Phi^*\Phi=\alpha I$ for some constant $\alpha > 0$. Define $g(x) = \pen(\Phi x)$, with $\pen$ being a lower semi-continuous, convex function. For any $x\in\CC^M$,
\begin{equation}
    \prox_{g}(x) = x + \frac{1}{\gamma}\Phi^*\left(\prox_{\gamma\pen}\left(\Phi x\right) - \Phi x \right).
    \label{eq:proxAx}
\end{equation}
\end{lemma}

\begin{figure}
		\definecolor{mycolor1}{rgb}{0.00000,0.44700,0.74100}%
		\definecolor{mycolor2}{rgb}{0.85000,0.32500,0.09800}%
		\definecolor{mycolor3}{rgb}{0.92900,0.69400,0.12500}%
		\definecolor{mycolor4}{rgb}{0.49400,0.18400,0.55600}%
		\definecolor{mycolor5}{rgb}{0,0,0}%
    \centering
		\footnotesize
		\begin{tikzpicture}
					\begin{axis}[
												width = 3in,
												height= 2in,
												xmin = 0,
												xmax = 5,
												xtick={0,2,4},
												xticklabels={0,$\lambda$,2$\lambda$},
												ytick={0,2,4},
												yticklabels={0,$\lambda$,2$\lambda$},
												samples=1600,
												ylabel near ticks,
												legend style={legend cell align=left, align=left, draw=black,at={(axis cs:0.1,5.25)},anchor=north west,font=\scriptsize,line width=0.2pt}
											]
							
							\addplot[solid,mycolor1,line width=1pt] {x*max(1-2/abs(x),0)};
							\addlegendentry{soft-thresholding}
							\addplot[solid,mycolor2,line width=1pt] {x*max(1-4/abs(x)^2,0)};
							\addlegendentry{non-neg.~Garrote}
							\addplot[solid,mycolor3,line width=1pt] {x*(abs(x)>2)};
							\addlegendentry{hard-thresholding}
							\draw [dashed, line width=1pt] ({axis cs:2,-0.2}|-{rel axis cs:0,0}) -- ({axis cs:2,5.2}|-{rel axis cs:0,5.2});
					\end{axis}
    \end{tikzpicture}
    \caption{Comparison of different thresholding rules for the identity function where $\lambda$ is the threshold.}
    \label{fig:nng}
\end{figure}
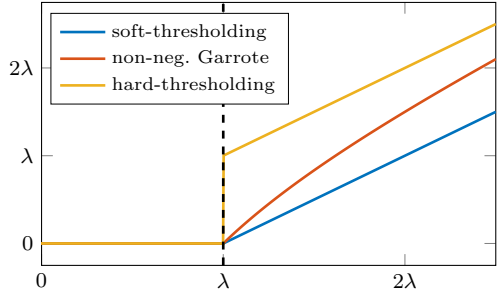

The following proximal operators are important in the remainder of this manuscript. First, with 
\begin{equation}
\pen_{\ell_1}(x) = \lambda\norm{x}_1 = \sum_{i=1}^N{\abs{x_i}},
\label{eq:l1norm}
\end{equation}
for some constant $\lambda>0$, the resulting proximity operator is the soft-thresholding operator
\begin{equation}
\prox_{\pen_{\ell_1}}(x) = x\cdot\max\left(1-\frac{\lambda}{\abs{x}},0\right).
\label{eq:st}
\end{equation}
Second, with $\lambda > 0$ setting
\begin{equation} 
\pen_\nng(x)=\sum_{i=1}^N\lambda^2 + \mbox{asinh}\left(\frac{\abs{x_i}}{2\lambda}\right)  + \lambda^2\frac{\abs{x_i}}{\sqrt{\abs{x_i}^2+4\lambda^2} + \abs{x_i}},
\label{eq:pennng}
\end{equation}
leads to the so called non-negative Garrote (NNG) \cite{Kowalski2014} defined by
\begin{equation}
\prox_{\pen_\nng}(x) = x\cdot\max\left(1-\frac{\lambda^2}{\abs{x}^2},0\right).
\label{eq:ew}
\end{equation}
Note that this proximal mapping is very similar to the soft-thresholding. It enforces sparsity and can be seen as a compromise between hard- and soft-thresholding, whereas the hard-thresholding for a given threshold $\lambda>0$ is defined by
\begin{equation}
    \mathcal{S}_h (x) = \left\{ \begin{array}{cc}
        x & \mbox{for~}\abs{x} > \lambda \\
        0 & \mbox{for~}\abs{x} \leq \lambda
    \end{array} \right. .
\end{equation}
This is visualized in Fig.~\ref{fig:nng} for the identity function. Function values after the non-negative Garrote are closer to the original values, but do not exhibit discontinuities as in the case of hard-thresholding.

\subsection{Undecimated wavelet transform}\label{sec:swt}
The wavelet transform is an orthogonal multi-scale transform that allows, in contrast to the ordinary Fourier transform, to resolve the spatial (or time) dependency of frequencies encoded in a signal. It allows to compress a large class of signals and is therefore a popular choice as a sparsity transformation in compressed-sensing applications \cite{lustig2007sparse}.

A single level of the discrete wavelet transform (DWT) for a one-dimensional signal $x$ consists of detail and approximation coefficients, $d$ and $a$ respectively. The detail and approximation coefficients result by filtering $x$ with low- and high-pass filter \textit{Lo} and \textit{Hi} followed by a dyadic downsampling \cite{mallat89}. A multilevel decomposition at level $j$ can then be realized by splitting the approximation coefficients $a^j$ again into corresponding detail and approximation coefficients $d^{j+1}$ and $a^{j+1}$. In a similar fashion the original signal can be reconstructed from its multilevel decomposition by a dyadic upsampling with zeros and convolving the results with the low- and high-pass reconstruction filters.

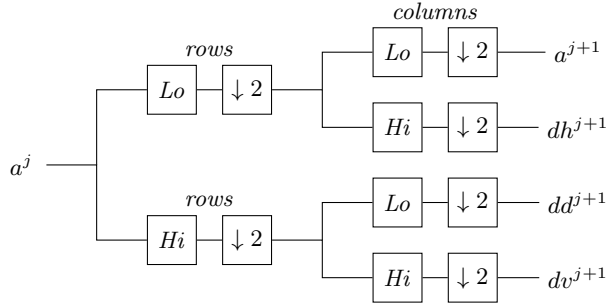
\begin{figure}
    \centering
    \small
    \begin{tikzpicture}
        \tikzset{every label/.style={xshift=5mm}}
        \tikzstyle{dd} = [coordinate]
        \tikzstyle{block} = [draw, rectangle, minimum height=2em, minimum width=2em]
        \tikzstyle{block2} = [rectangle, minimum height=2em, minimum width=2em]
        
        \node[block2,name=A] at (0,2.5){$a^j$};
        \node[dd,name=A1] at (1,2.5){};
        \node[dd,name=A2] at (1,1.5){};
        \node[dd,name=A3] at (1,3.5){};
        \node[block,name=B,label=\small\emph{rows}] at (2,1.5){\textit{Hi}};
        \node[block,name=C,label=\small\emph{rows}] at (2,3.5){\textit{Lo}};
        \node[block,name=D] at (3,1.5){$\downarrow2$};
        \node[dd,name=B1] at (4,1.5){};
        \node[dd,name=B2] at (4,1){};
        \node[dd,name=B3] at (4,2){};
        \node[block,name=E] at (3,3.5){$\downarrow2$};
        \node[dd,name=C1] at (4,3.5){};
        \node[dd,name=C2] at (4,3){};
        \node[dd,name=C3] at (4,4){};
        \node[block,name=F,label=\small\emph{columns}] at (5,4){\textit{Lo}};
        \node[block,name=G] at (5,3){\textit{Hi}};
        \node[block,name=H] at (5,2){\textit{Lo}};
        \node[block,name=I] at (5,1){\textit{Hi}};
        \node[block,name=J] at (6,4){$\downarrow2$};
        \node[block,name=K] at (6,3){$\downarrow2$};        
        \node[block,name=L] at (6,1){$\downarrow2$};
        \node[block,name=M] at (6,2){$\downarrow2$};
        \node[block2,name=N] at (7.4,4){$a^{j+1}$};
        \node[block2,name=O] at (7.4,3){$dh^{j+1}$};
        \node[block2,name=P] at (7.4,1){$dv^{j+1}$};
        \node[block2,name=Q] at (7.4,2){$dd^{j+1}$};
        
        \draw[draw,-] (A)--(A1)--(A2)--(B)--(D)--(B1)--(B2)--(I)--(L)--(P);
        \draw[draw,-] (A1)--(A3)--(C)--(E)--(C1)--(C2)--(G)--(K)--(O);
        \draw[draw,-] (B1)--(B3)--(H)--(M)--(Q);
        \draw[draw,-] (C1)--(C3)--(F)--(J)--(N);

    \end{tikzpicture}
    \caption{Schematic of a single level 2D DWT decomposition. A block with $\downarrow 2$ denotes dyadic downsampling, whereas the others denote the convolution with low- and high-pass filter coefficients along rows or columns.}
    \label{fig:dwt}
\end{figure}

The DWT can be extended to higher dimensional signals by tensor products. For a two-dimensional signal, for example, a single level DWT decomposition results in four components: the approximation coefficients $a^j$ and three detail coefficients $dh^j$, $dv^j$ and $dd^d$ corresponding to vertical, horizontal and diagonal details. The schematic principle is visualized in Fig.~\ref{fig:dwt}.

In order to reconstruct the original signal from its wavelet coefficients, the steps in Fig.~\ref{fig:dwt} are simply reversed. That means each of the four components are upsampled, filtered with the appropriate reconstruction filter and summed up \cite{mallat89}. In the following we denote the reconstruction at level $j$ by
\begin{equation}
    a^j = \mathcal{R}^j\left(a^{j+1},dh^{j+1},dv^{j+1},dv^{j+1}\right).
\end{equation}

Although the DWT is applied in numerous image processing applications, e.g., compressed sensing in magnetic resonance imaging, it suffers from a lack of shift-invariance and aliasing that results from the inherent downsampling. Omitting the downsampling step results in the undecimated discrete wavelet transform (UDWT), which is also named à trous or stationary wavelet transform. Instead of downsampling, filter coefficients are upsampled by inserting zeros at even- or odd-indexed positions. This introduces a redundancy of the wavelet coefficients making reconstruction a bit more complex. It basically involves reconstructing the even and odd parts of the components separately followed by an averaging step, which combines both reconstructed parts. For simplifying notation we define $c_j\coloneqq\left\{a^{j+1},dh^{j+1},dv^{j+1},dv^{j+1}\right\}$. Furthermore, let the operators that select every even and odd member of a doubly indexed sequence $x_{m,n}$ be $(D_0x)_{m,n} = x_{2m,2n}$ and $(D_1 x)_{m,n} = x_{2m+1,2n+1}$. The original data then can be reconstructed by recursively computing
\begin{equation}
    c^j = \frac{1}{2}\left(\mathcal{R}^j( D_0c^{j+1}) + \mathcal{R}^j(D_1c^{j+1}) \right).
    \label{eq:iudwt}
\end{equation}
For more details we refer to the concept of the $\epsilon$-decimated DWT in \cite{Nason1995}. 

The extension of the 2D UDWT to three dimensions results in eight components, by low- and highpass filtering each of the four 2D components along the third dimension. The scheme in \eqref{eq:iudwt} can be extended in a similar fashion to reconstruct the original 3D signal.

\subsection{Kaczmarz Algorithm}
In the following let $A\in\CC^{m\times n}$ be a matrix with rows $a^T_i\in\CC^n$, $i=1,2,\ldots,m$ and $b\in\CC^m$.  
The standard Kaczmarz algorithm, initially proposed in \cite{Kaczmarz37}, solves the consistent linear system of equations $Ax=b$, such that the solution $x^*$ satisfies 
$$
x^* = \argmin_x \norm{x}_2\quad \mbox{s.t.~} Ax=b.
$$ 
Given an initial $x_0$, the Kaczmarz algorithm reads
\begin{equation}
x_{k+1}= x_{k}+\frac{b_{i} - \langle a_{i} \,,\, x_{k}\rangle}{\Vert a_{i}\Vert _{2}^{2}}a_i^*, \quad\quad i=\mbox{mod}(k,m).
\label{eq:kz}
\end{equation}
If $x_0=0$, the Kaczmarz algorithm converges to a minimum norm solution \cite{Schoepfer2018}.
The rate of convergence depends on the ordering of the rows and is shown to be exponential in a randomized Kaczmarz variant, where the rows are chosen randomly \cite{Strohmer2008}.
As it only uses a single row in each iteration, it is called a \textit{row-action method}. 
Instead of considering the constraint $Ax=b$ as an intersection of several smaller linear constraints, the iteration in \eqref{eq:kz} can be directly related to the Landweber method
\begin{equation}
x_{k+1} = x_k + \gamma A^*\left(Ax_k - b\right),
\label{eq:landweber}
\end{equation}
with $\gamma=\norm{A}^{-2}_2$, see for example \cite[Sec.~2.5.1]{Lorenz2014a} or \cite{Trussell1985}.

\section{Proposed Algorithm}\label{sec:propAlg}

As already mentioned in \cite{Storath2017}, an $\ell_2$-prior which penalizes the total energy of the tracer distribution does not take the structure of the latent image into account. Furthermore, this reconstruction model assumes a purely Gaussian noise pattern independently of each pixel intensity. In the case of MPI this is simply not a very accurate assumption. First, there are always image parts with no tracer concentration present and hence all these background pixels are correlated. On the other hand, an equally distributed tracer concentration in the sample, i.e., a state of equilibrium, leads to images with piecewise constant regions. Both assumptions apparently motivate the usage of an $\ell_1$-norm and total variation norm regularization approach \cite{Storath2017, ilbey2017}.

Similarly, piecewise regular signals are also sparsely represented in a discrete wavelet basis. As these wavelets are well localized, only few coefficients are needed to represent transient structures such as edges. Hence, tracer concentrations in Magnetic Particle Imaging exhibit large wavelet coefficients in the vicinity of large transients and small wavelet coefficients for regular textures over the support of the corresponding wavelet \cite{mallat89}. This motivates the following optimization problem.    

\begin{algorithm}[!t]
\setstretch{1.3}
\SetAlgoLined
\DontPrintSemicolon
\SetKwInOut{Input}{Input}
\SetKwInOut{Output}{Output}
\vspace{0.1cm}
\Input{$A \in\CC^{m\times n}$ - measured system matrix \\ 
       $b \in \CC^m$ - measured phantom signal\\
			 $\lambda>0$ - regularization parameter\\
			 $\gamma>0$ - Lipschitz constant of $A$\\
			 $\Phi$ - tight wavelet frame\\
			 $N$ - maximum iteration number }
\Output{$x \in \RR^{n}$ - reconstructed tracer concentration}		
\BlankLine
	Fix $x_0 \in \CC^{n}$, $z_0 \gets x_0$ and $t_0 = 1$\;
	\For{$k=1,2,\ldots,N$}
	{
		$y_k \gets  z_{k-1} - \frac{1}{\gamma} A^*(Az_{k-1} - b)$\;
		$y_k \gets P(y_k)$\tcp*{\footnotesize Projection onto $\RR^+$}
		$x_k \gets \prox_{\gamma^{-1}\lambda\pen}(\Phi y_{k})$\tcp*{\footnotesize see Lemma \ref{lem:proxA}}
		$t_k \gets \frac{1 + \sqrt{4t_{k-1}^2+1}}{2}$\;
		$z_k \gets x_{k-1} + \left(1+\frac{t_{k-1} - 1}{t_k} \right)\left(x_k - x_{k-1}\right)$\;
	}
	\vspace{0.1cm}
\caption{FISTA}
\label{alg:fista}
\end{algorithm}

With $\Phi:\RR^n \rightarrow \RR^{J\times n}$, denoting the undecimated wavelet transform described above, with $J$ being the number of decomposition levels, we propose the following minimization problem for reconstructing the tracer concentration in MPI:
\begin{equation}
\min_{x\in\RR^+} \frac{1}{2}\norm{Ax-b}_2^2 + \lambda \pen(\Phi x).
\label{eq:minpr}
\end{equation}
Here, $\pen$ promotes the sparsity of the wavelet coefficients $\Phi x$, either via the simple $\ell_1$-norm from \eqref{eq:l1norm} or the NNG representation from \eqref{eq:ew}. 

There exist several algorithms for solving problem \eqref{eq:minpr}, such as FISTA \cite{Beck2009a}, ADMM \cite{Boyd2010} or Douglas-Rachford-splitting approaches \cite{Combettes2011}. 
The FISTA approach is based on a general forward-backward splitting approach by iterating
\begin{equation}
x_{k+1} = \prox_{\gamma \pen}\left( x_{k} - \gamma \nabla h(x_k) \right),\qquad\gamma\in\ ]0,\infty[\ ,
\label{eq:fbs}
\end{equation}
with $h(x) = \frac{1}{2}\norm{Ax-b}^2_2$. This scheme is characterized by its gradient (forward) step and its backward step based on the proximity operator of $\pen$. The resulting FISTA based algorithm is summarized in Algorithm \ref{alg:fista}. Note that in practice $\gamma = \norm{A^*A}_{\text{op}}^{-1}=\frac{1}{\sqrt{\abs{\varrho}}}$, where $\varrho$ is the largest eigenvalue of $A^*A$.

Our proposed approach for solving \eqref{eq:minpr} is described in the following. Based on the forward-backward splitting in \eqref{eq:fbs}, we propose to use the Kaczmarz algorithm for the gradient step. To be more specific, we propose a single sweep over all rows of the system matrix $A$ to compute a gradient descent step. This is justified by the direct relation of the Landweber method in \eqref{eq:landweber} and the Kaczmarz iteration \eqref{eq:kz}. After that we simply apply the proximal operator corresponding to $\pen$ resulting in the following iteration  
\begin{equation}
\begin{aligned} 
		x_{k} &= x_{k}+\frac{b_{i} - \langle a_{i} \,,\, x_{k}\rangle}{\Vert a_{i}\Vert _{2}^{2}} a^*_{i}, \quad i=\mbox{mod}(k,m) \\[8pt]
		x_{k+1}&= \left\{ \begin{array}{ll}  \prox_{\pen}(\Phi x_{k})& \mbox{if }i = 0\\ x_{k} & \mbox{if }i \neq 0\end{array}  \right.
	\end{aligned} 
	\label{eq:3}
\end{equation}
which leads to the following pseudocode of our proposed Sparse Kaczmarz Algorithm (SKA) summarized in Algorithm \ref{alg:sKA}. Note that the tracer concentration in MPI is nonnegative implying to iteratively project onto the positive (real) numbers \cite{Kluth2019}.

\begin{algorithm}[!t]
\setstretch{1.3}
\SetAlgoLined
\DontPrintSemicolon
\SetKwInOut{Input}{Input}
\SetKwInOut{Output}{Output}
\vspace{0.1cm}
\Input{$A \in\CC^{m\times n}$ - measured system matrix \\ 
       $b \in \CC^m$ - measured phantom signal\\
			 $\lambda>0$ - regularization parameter\\
			 $\Phi$ - tight wavelet frame\\
			 $N$ - maximum iteration number }
\Output{$x \in \RR^{n}$ - reconstructed tracer concentration}		
\BlankLine
	$x_0 \gets \mathbf{0}\in\CC^{n}$\;
	\For{$k=1,2,\ldots,N$}
	{
		\For{$i=1,2,\ldots,m$}
		{
		  $x_{k} \gets x_{k}+\frac{b_{i} - \langle a_{i} \,,\, x_{k}\rangle }{\Vert a_{i}\Vert _{2}^{2} }a_{i}^*$\;
	  }
		$x_k \gets P(y_x)$\tcp*{\footnotesize Projection onto $\RR^+$}
		$x_k \gets \mbox{prox}_{\lambda\pen}\left(\Phi x_k\right)$\tcp*{\footnotesize see Lemma \ref{lem:proxA}}
	}
	\vspace{0.1cm}
\caption{Sparse Kaczmarz Algorithm (SKA)}
\label{alg:sKA}
\end{algorithm}

There are several issues that need clarification at this point. First of all, the Kaczmarz algorithm in general only works for consistent linear system, i.e., it must hold that $b\in\mathcal{R}(A)$ where $\mathcal{R}$ denotes the range of $A$. Obviously, in the case of MPI it holds that $b\notin\mathcal{R}(A)$ due to measurement errors as well as noise. In such cases, convergence of the Kaczmarz iteration for inconsistent linear systems is proven in \cite[Thm.~3.2]{Chen2018}, with the rate of convergence depending on the amount of noise. In MPI, the amount of noise present can be, at least to some extent, controlled by preconditioning the system matrix via SNR thresholding \cite{Knopp2010} as described above. This is similar to the so called bang-bang relaxation parameter proposed in \cite{Scherzer2007} and justifies the usage of a basic Kaczmarz iteration. 

Another issue is the convergence of $\pen_\nng$ in the context of general forward-backward splitting approaches as it is not strictly convex. However, the convergence of general semi-convex (weakly-convex) penalty functions inside an iterative shrinkage/thresholding algorithm (ISTA) has been shown in \cite{Bayram2016,Kowalski2014}, inside an Douglas-Rachford splitting in \cite{Guo2017} and inside the alternating direction method of multipliers (ADMM) recently in \cite{Zhang2019}. 

In order to directly compare the first order proximal method in Algorithm \ref{alg:fista} with the Kaczmarz based Algorithm \ref{alg:sKA}, the Kaczmarz based gradient descend has to adopt a single step size for all the rows of the system matrix. This requires a row normalization of $A$ such that $\norm{a_i}_2 = 1$ for all $i\in\{1,2,\ldots,m\}$.

\section{Numerical Results} \label{sec:res}
\subsection{Experimental Setup}

MPI data is simulated using a measured hybrid system matrix \cite{Gladiss2017, Graeser2017}. The excitation field amplitude is set to 12\,mT in both $x$- and $y$-directions with a frequency of 24.5098\,kHz in $x$- and 26.0417\,kHz in $y$-direction at a discretization of 0.5\,mT. The FOV is discretized at 57$\times$57 sampling points. The bandwidth of the receiver unit is 4.375\,MHz. The induced signals in the $x$- and $y$-receive coils are averaged over 50 repetitions, Fourier transformed and stacked as the columns of the system matrix $A$. Furthermore, $A$ is row-normalized such that all rows have unit norm. The resulting system matrix has a size of 5714$\times$3249. 
Forward simulations are based on 
\begin{equation}
b = A(\sigma^{-1} x )+ \eta,
\label{eq:fwdsimul}
\end{equation}
with a predefined phantom $x$ with weight $\sigma$ and a noise vector $\eta$ consisting of the mean over several empty scanner measurements. 
In contrast to the simulations in \cite{Storath2017, Bathke2017}, where only white Gaussian noise was added, this approach is much closer to real data, since the noise in an MPI experiment is typically colored due to a non-constant transfer function of the receiver electronics. In our case, the actual noise level is defined by $\sigma$, which implies that the larger $\sigma$, the larger the influence of the constant background signal. Furthermore, from the resulting simulated signal $b$ only frequencies within the range of 70-3,000\,kHz are considered in the reconstruction process \cite{Knopp2010,Kluth2019}. This reduces the number of frequency components from 5712 to 3065 and results in a total system matrix size of 3065$\times$3249.

In the following forward simulations we consider the phantoms visualized in Fig.~\ref{fig:resC}. The generic shape phantom (Fig.~\ref{fig:resC}a) is taken from \cite{Hansen1994} and has been used in the context of MPI in \cite{Bathke2017}. The vascular tree phantom shown in Fig.~\ref{fig:resC}b has a more medically motivated structure, driven by angiographical applications. The quality of the reconstructed tracer concentration $x^{\text{rec}}$ is evaluated using the peak signal-to-noise ratio (PSNR) 
\begin{equation}
\mbox{PSNR} (x^{\text{rec}}) = 10 \log_{10}\left( \frac{1}{n}\sum\limits_{i=1}^{n}\left(\sigma x^{\text{rec}}_i - x_i \right)^2\right)^{-1},
\label{eq:psnr}
\end{equation}
and the structural similarity index SSIM \cite{Wang2004}, where SSIM values closer to 1.0 represents a better perceived image quality between original and reconstructed image.

The maximum iteration number of all four algorithms is set to 3,000 with the following stopping criterion: whenever the relative change between two consecutive iterates $\epsrel$ is below $10^{-5}$, i.e., 
\begin{equation}
\epsrel=\frac{\norm{x_k-x_{k+1}}_2}{\norm{x_k}_2} < 10^{-5},
\label{eq:conv}
\end{equation}
the iteration is stopped \cite{Storath2017,ilbey2017}. The regularization parameter $\lambda$ is chosen based on a greedy approach maximizing the corresponding PSNR. The underlying undecimated wavelet transform is based on the Haar wavelet and a decomposition level of two.
In the following we will consider the phantom with a low noise level $\sigma=10$ and a high noise level $\sigma=50$.
Computation times are based on a workstation with 32GB RAM and a 2.9GHz i7 QuadCore CPU.

\definecolor{mycolor1}{rgb}{0.00000,0.44700,0.74100}%
\definecolor{mycolor2}{rgb}{0.85000,0.32500,0.09800}%
\definecolor{mycolor3}{rgb}{0.92900,0.69400,0.12500}%
\definecolor{mycolor4}{rgb}{0.49400,0.18400,0.55600}%

\newcommand{\myPSNRplot}[4]{
\scriptsize
\begin{tikzpicture}
				\begin{axis}[%
				width=2in,
				height=1.5in,
				at={(0.0in,0.0in)},
				scale only axis,
				separate axis lines,
				every outer x axis line/.append style={black},
				every x tick label/.append style={font=\color{black}},
				every x tick/.append style={black},
				xmin=0,
				xmax=#1,
				xlabel near ticks,
				xlabel={Computation time (s)},
				every outer y axis line/.append style={black},
				every y tick label/.append style={font=\color{black}},
				every y tick/.append style={black},
				ymin=#2,
				ymax=#3,
				ylabel={PSNR},
				ylabel near ticks,
				axis background/.style={fill=white},
				legend style={at={(0.97,0.03)}, anchor=south east, legend cell align=left, align=left, draw=black,font={\notsotiny},line width = 0.2pt}
				]
				\addplot [color=mycolor1, line width=0.7pt]
					table[]{figures/data/resPSNR_#4-1.tsv};
				\addlegendentry{SKA (NNG)}

				\addplot [color=mycolor2, line width=0.7pt]
					table[]{figures/data/resPSNR_#4-2.tsv};
				\addlegendentry{FISTA (NNG)}

				\addplot [dashed, color=mycolor1, line width=0.7pt]
					table[]{figures/data/resPSNR_#4-3.tsv};
				\addlegendentry{SKA (ST)}

				\addplot [dashed, color=mycolor2, line width=0.7pt]
					table[]{figures/data/resPSNR_#4-4.tsv};
				\addlegendentry{FISTA (ST)}

				\end{axis}
		\end{tikzpicture}
}

\begin{figure}[t!]
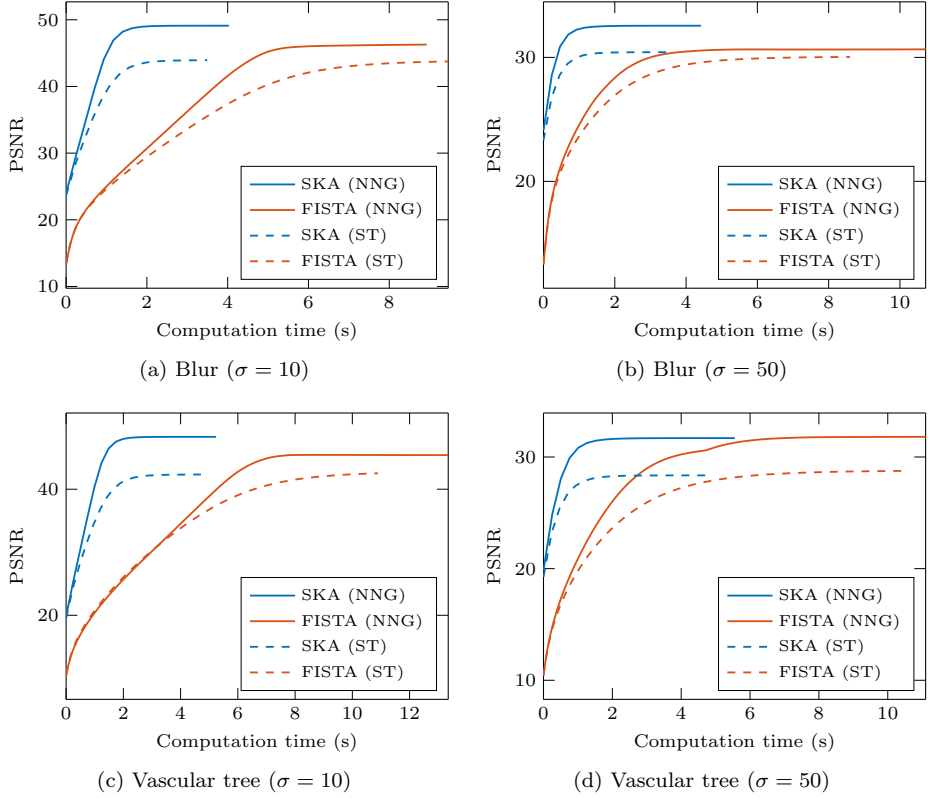
%
    \centering
	\subfloat[Blur ($\sigma=10$)]{
		\myPSNRplot{9.46}{9.74}{52.69}{blur01}
	\label{fig:res1_1}
	}
	\subfloat[Blur ($\sigma=50$)]{
		\myPSNRplot{10.72}{11.38}{34.48}{blur002}
	\label{fig:res1_2}
	}\\

	\subfloat[Vascular tree ($\sigma=10$)]{
		\myPSNRplot{13.35}{6.66}{52.09}{vasc01}
	\label{fig:res1_3}
	}
	\subfloat[Vascular tree ($\sigma=50$)]{
		\myPSNRplot{11.1}{8.29}{33.94}{vasc002}
	\label{fig:res1_4}
	}
\caption{PSNR values versus the computation time for the shape and vascular tree phantom at low noise level $\sigma=10$ and large noise level $\sigma=50$.}%
\label{fig:res1}%
\end{figure}

\begin{table*}[t!]
\renewcommand{\arraystretch}{1.3}
\caption{Quality measures (psnr, ssim and computation time in seconds) of the reconstructed particle concentrations based on the proposed sparse Kaczmarz variants (SKA), the FISTA approaches and current state-of-the-art MPI reconstruction algorithms for the shape and vascular tree phantom at different noise levels $\sigma$. }
\label{tab:1}
\footnotesize
\vspace{10pt}
\begin{tabular}{@{}lllllllllllll@{}} \toprule
& \multicolumn{3}{c}{Shape ($\sigma=10$)} & \multicolumn{3}{c}{Vascular tree ($\sigma=10$)} & \multicolumn{3}{c}{Shape ($\sigma=50$)} & \multicolumn{3}{c}{Vascular tree ($\sigma=50$)} \\ \cmidrule(r){2-4} \cmidrule(r){5-7} \cmidrule(r){8-10} \cmidrule(r){11-13}
															 & PSNR  & SSIM   & time & PSNR  & SSIM   & time & PSNR  & SSIM   & time & PSNR  & SSIM   & time  \\ \midrule
SKA (NNG)   									 & 49.11 & 0.9887 & 4.0  & 48.31 & 0.9955 & 5.2  & 32.56 & 0.8256 & 4.4  & 31.69 & 0.8873 & 5.6   \\
FISTA (NNG) 									 & 46.28 & 0.9115 & 8.9  & 45.41 & 0.9471 & 13.3 & 30.65 & 0.6991 & 8.9  & 31.80 & 0.7840 & 11.1  \\ \midrule
SKA (ST)    									 & 43.93 & 0.9537 & 3.5  & 42.33 & 0.9525 & 4.7  & 30.42 & 0.7315 & 3.6  & 28.35 & 0.7268 & 4.9   \\
FISTA (ST)                     & 43.75 & 0.8871 & 9.5  & 42.52 & 0.9220 & 10.9 & 30.04 & 0.6806 & 8.6  & 28.76 & 0.7791 & 10.6  \\\midrule
reg. KA \cite{Knopp2010}       & 38.49 & 0.8475 & 122  & 39.55 & 0.8899 & 135  & 25.22 & 0.4767 & 195  & 25.27 & 0.6526 & 189   \\
Fused Lasso \cite{Storath2017} & 46.02 & 0.9421 & 304  & 44.46 & 0.9463 & 310  & 29.81 & 0.8614 & 324  & 30.91 & 0.8955 & 321  \\\bottomrule
\end{tabular}
\end{table*}

\newcommand{\myResPlot}[2]{
\scriptsize
\begin{tikzpicture}
\begin{axis}[%
width=1in,
height=1in,
at={(0in,0in)},
scale only axis,
point meta min=0,
point meta max=#1,
axis on top,
separate axis lines,
every outer x axis line/.append style={black},
every x tick label/.append style={font=\color{black}},
every x tick/.append style={black},
xmin=0.5,
xmax=57.5,
xlabel = {Original},
xlabel near ticks,
every outer y axis line/.append style={black},
every y tick label/.append style={font=\color{black}},
every y tick/.append style={black},
y dir=reverse,
ymin=0.5,
ymax=57.5,
axis line style={draw=none},
ticks=none,
legend style={legend cell align=left, align=left, draw=black},
]
\addplot [forget plot] graphics [xmin=0.5, xmax=57.5, ymin=0.5, ymax=57.5] {figures/#2-1.png};
\end{axis}

\begin{axis}[%
width=1in,
height=1in,
at={(1.075in,0in)},
scale only axis,
point meta min=0,
point meta max=#1,
axis on top,
separate axis lines,
every outer x axis line/.append style={black},
every x tick label/.append style={font=\color{black}},
every x tick/.append style={black},
xmin=0.5,
xmax=57.5,
xlabel = {SKA (ST)},
xlabel near ticks,
every outer y axis line/.append style={black},
every y tick label/.append style={font=\color{black}},
every y tick/.append style={black},
y dir=reverse,
ymin=0.5,
ymax=57.5,
axis line style={draw=none},
ticks=none,
legend style={legend cell align=left, align=left, draw=black},
]
\addplot [forget plot] graphics [xmin=0.5, xmax=57.5, ymin=0.5, ymax=57.5] {figures/#2-4.png};
\end{axis}

\begin{axis}[%
width=1in,
height=1in,
at={(2.15in,0in)},
scale only axis,
point meta min=0,
point meta max=#1,
axis on top,
separate axis lines,
every outer x axis line/.append style={black},
every x tick label/.append style={font=\color{black}},
every x tick/.append style={black},
xmin=0.5,
xmax=57.5,
xlabel = {FISTA (ST)},
xlabel near ticks,
every outer y axis line/.append style={black},
every y tick label/.append style={font=\color{black}},
every y tick/.append style={black},
y dir=reverse,
ymin=0.5,
ymax=57.5,
axis line style={draw=none},
ticks=none,
legend style={legend cell align=left, align=left, draw=black},
]
\addplot [forget plot] graphics [xmin=0.5, xmax=57.5, ymin=0.5, ymax=57.5] {figures/#2-5.png};
\end{axis}

\begin{axis}[%
width=1in,
height=1in,
at={(3.225in,0in)},
scale only axis,
point meta min=0,
point meta max=#1,
axis on top,
separate axis lines,
every outer x axis line/.append style={black},
every x tick label/.append style={font=\color{black}},
every x tick/.append style={black},
xmin=0.5,
xmax=57.5,
xlabel = {FISTA (NNG)},
xlabel near ticks,
every outer y axis line/.append style={black},
every y tick label/.append style={font=\color{black}},
every y tick/.append style={black},
y dir=reverse,
ymin=0.5,
ymax=57.5,
axis line style={draw=none},
ticks=none,
legend style={legend cell align=left, align=left, draw=black},
]
\addplot [forget plot] graphics [xmin=0.5, xmax=57.5, ymin=0.5, ymax=57.5] {figures/#2-3.png};
\end{axis}

\begin{axis}[%
width=1in,
height=1in,
at={(4.3in,0in)},
scale only axis,
point meta min=0,
point meta max=#1,
axis on top,
separate axis lines,
every outer x axis line/.append style={black},
every x tick label/.append style={font=\color{black}},
every x tick/.append style={black},
xmin=0.5,
xmax=57.5,
xlabel = {SKA (NNG)},
xlabel near ticks,
every outer y axis line/.append style={black},
every y tick label/.append style={font=\color{black}},
every y tick/.append style={black},
y dir=reverse,
ymin=0.5,
ymax=57.5,
axis line style={draw=none},
ticks=none,
legend style={legend cell align=left, align=left, draw=black},
colormap={mymap}{[1pt] rgb(0pt)=(0.2422,0.1504,0.6603); rgb(1pt)=(0.2444,0.1534,0.6728); rgb(2pt)=(0.2464,0.1569,0.6847); rgb(3pt)=(0.2484,0.1607,0.6961); rgb(4pt)=(0.2503,0.1648,0.7071); rgb(5pt)=(0.2522,0.1689,0.7179); rgb(6pt)=(0.254,0.1732,0.7286); rgb(7pt)=(0.2558,0.1773,0.7393); rgb(8pt)=(0.2576,0.1814,0.7501); rgb(9pt)=(0.2594,0.1854,0.761); rgb(11pt)=(0.2628,0.1932,0.7828); rgb(12pt)=(0.2645,0.1972,0.7937); rgb(13pt)=(0.2661,0.2011,0.8043); rgb(14pt)=(0.2676,0.2052,0.8148); rgb(15pt)=(0.2691,0.2094,0.8249); rgb(16pt)=(0.2704,0.2138,0.8346); rgb(17pt)=(0.2717,0.2184,0.8439); rgb(18pt)=(0.2729,0.2231,0.8528); rgb(19pt)=(0.274,0.228,0.8612); rgb(20pt)=(0.2749,0.233,0.8692); rgb(21pt)=(0.2758,0.2382,0.8767); rgb(22pt)=(0.2766,0.2435,0.884); rgb(23pt)=(0.2774,0.2489,0.8908); rgb(24pt)=(0.2781,0.2543,0.8973); rgb(25pt)=(0.2788,0.2598,0.9035); rgb(26pt)=(0.2794,0.2653,0.9094); rgb(27pt)=(0.2798,0.2708,0.915); rgb(28pt)=(0.2802,0.2764,0.9204); rgb(29pt)=(0.2806,0.2819,0.9255); rgb(30pt)=(0.2809,0.2875,0.9305); rgb(31pt)=(0.2811,0.293,0.9352); rgb(32pt)=(0.2813,0.2985,0.9397); rgb(33pt)=(0.2814,0.304,0.9441); rgb(34pt)=(0.2814,0.3095,0.9483); rgb(35pt)=(0.2813,0.315,0.9524); rgb(36pt)=(0.2811,0.3204,0.9563); rgb(37pt)=(0.2809,0.3259,0.96); rgb(38pt)=(0.2807,0.3313,0.9636); rgb(39pt)=(0.2803,0.3367,0.967); rgb(40pt)=(0.2798,0.3421,0.9702); rgb(41pt)=(0.2791,0.3475,0.9733); rgb(42pt)=(0.2784,0.3529,0.9763); rgb(43pt)=(0.2776,0.3583,0.9791); rgb(44pt)=(0.2766,0.3638,0.9817); rgb(45pt)=(0.2754,0.3693,0.984); rgb(46pt)=(0.2741,0.3748,0.9862); rgb(47pt)=(0.2726,0.3804,0.9881); rgb(48pt)=(0.271,0.386,0.9898); rgb(49pt)=(0.2691,0.3916,0.9912); rgb(50pt)=(0.267,0.3973,0.9924); rgb(51pt)=(0.2647,0.403,0.9935); rgb(52pt)=(0.2621,0.4088,0.9946); rgb(53pt)=(0.2591,0.4145,0.9955); rgb(54pt)=(0.2556,0.4203,0.9965); rgb(55pt)=(0.2517,0.4261,0.9974); rgb(56pt)=(0.2473,0.4319,0.9983); rgb(57pt)=(0.2424,0.4378,0.9991); rgb(58pt)=(0.2369,0.4437,0.9996); rgb(59pt)=(0.2311,0.4497,0.9995); rgb(60pt)=(0.225,0.4559,0.9985); rgb(61pt)=(0.2189,0.462,0.9968); rgb(62pt)=(0.2128,0.4682,0.9948); rgb(63pt)=(0.2066,0.4743,0.9926); rgb(64pt)=(0.2006,0.4803,0.9906); rgb(65pt)=(0.195,0.4861,0.9887); rgb(66pt)=(0.1903,0.4919,0.9867); rgb(67pt)=(0.1869,0.4975,0.9844); rgb(68pt)=(0.1847,0.503,0.9819); rgb(69pt)=(0.1831,0.5084,0.9793); rgb(70pt)=(0.1818,0.5138,0.9766); rgb(71pt)=(0.1806,0.5191,0.9738); rgb(72pt)=(0.1795,0.5244,0.9709); rgb(73pt)=(0.1785,0.5296,0.9677); rgb(74pt)=(0.1778,0.5349,0.9641); rgb(75pt)=(0.1773,0.5401,0.9602); rgb(76pt)=(0.1768,0.5452,0.956); rgb(77pt)=(0.1764,0.5504,0.9516); rgb(78pt)=(0.1755,0.5554,0.9473); rgb(79pt)=(0.174,0.5605,0.9432); rgb(80pt)=(0.1716,0.5655,0.9393); rgb(81pt)=(0.1686,0.5705,0.9357); rgb(82pt)=(0.1649,0.5755,0.9323); rgb(83pt)=(0.161,0.5805,0.9289); rgb(84pt)=(0.1573,0.5854,0.9254); rgb(85pt)=(0.154,0.5902,0.9218); rgb(86pt)=(0.1513,0.595,0.9182); rgb(87pt)=(0.1492,0.5997,0.9147); rgb(88pt)=(0.1475,0.6043,0.9113); rgb(89pt)=(0.1461,0.6089,0.908); rgb(90pt)=(0.1446,0.6135,0.905); rgb(91pt)=(0.1429,0.618,0.9022); rgb(92pt)=(0.1408,0.6226,0.8998); rgb(93pt)=(0.1383,0.6272,0.8975); rgb(94pt)=(0.1354,0.6317,0.8953); rgb(95pt)=(0.1321,0.6363,0.8932); rgb(96pt)=(0.1288,0.6408,0.891); rgb(97pt)=(0.1253,0.6453,0.8887); rgb(98pt)=(0.1219,0.6497,0.8862); rgb(99pt)=(0.1185,0.6541,0.8834); rgb(100pt)=(0.1152,0.6584,0.8804); rgb(101pt)=(0.1119,0.6627,0.877); rgb(102pt)=(0.1085,0.6669,0.8734); rgb(103pt)=(0.1048,0.671,0.8695); rgb(104pt)=(0.1009,0.675,0.8653); rgb(105pt)=(0.0964,0.6789,0.8609); rgb(106pt)=(0.0914,0.6828,0.8562); rgb(107pt)=(0.0855,0.6865,0.8513); rgb(108pt)=(0.0789,0.6902,0.8462); rgb(109pt)=(0.0713,0.6938,0.8409); rgb(110pt)=(0.0628,0.6972,0.8355); rgb(111pt)=(0.0535,0.7006,0.8299); rgb(112pt)=(0.0433,0.7039,0.8242); rgb(113pt)=(0.0328,0.7071,0.8183); rgb(114pt)=(0.0234,0.7103,0.8124); rgb(115pt)=(0.0155,0.7133,0.8064); rgb(116pt)=(0.0091,0.7163,0.8003); rgb(117pt)=(0.0046,0.7192,0.7941); rgb(118pt)=(0.0019,0.722,0.7878); rgb(119pt)=(0.0009,0.7248,0.7815); rgb(120pt)=(0.0018,0.7275,0.7752); rgb(121pt)=(0.0046,0.7301,0.7688); rgb(122pt)=(0.0094,0.7327,0.7623); rgb(123pt)=(0.0162,0.7352,0.7558); rgb(124pt)=(0.0253,0.7376,0.7492); rgb(125pt)=(0.0369,0.74,0.7426); rgb(126pt)=(0.0504,0.7423,0.7359); rgb(127pt)=(0.0638,0.7446,0.7292); rgb(128pt)=(0.077,0.7468,0.7224); rgb(129pt)=(0.0899,0.7489,0.7156); rgb(130pt)=(0.1023,0.751,0.7088); rgb(131pt)=(0.1141,0.7531,0.7019); rgb(132pt)=(0.1252,0.7552,0.695); rgb(133pt)=(0.1354,0.7572,0.6881); rgb(134pt)=(0.1448,0.7593,0.6812); rgb(135pt)=(0.1532,0.7614,0.6741); rgb(136pt)=(0.1609,0.7635,0.6671); rgb(137pt)=(0.1678,0.7656,0.6599); rgb(138pt)=(0.1741,0.7678,0.6527); rgb(139pt)=(0.1799,0.7699,0.6454); rgb(140pt)=(0.1853,0.7721,0.6379); rgb(141pt)=(0.1905,0.7743,0.6303); rgb(142pt)=(0.1954,0.7765,0.6225); rgb(143pt)=(0.2003,0.7787,0.6146); rgb(144pt)=(0.2061,0.7808,0.6065); rgb(145pt)=(0.2118,0.7828,0.5983); rgb(146pt)=(0.2178,0.7849,0.5899); rgb(147pt)=(0.2244,0.7869,0.5813); rgb(148pt)=(0.2318,0.7887,0.5725); rgb(149pt)=(0.2401,0.7905,0.5636); rgb(150pt)=(0.2491,0.7922,0.5546); rgb(151pt)=(0.2589,0.7937,0.5454); rgb(152pt)=(0.2695,0.7951,0.536); rgb(153pt)=(0.2809,0.7964,0.5266); rgb(154pt)=(0.2929,0.7975,0.517); rgb(155pt)=(0.3052,0.7985,0.5074); rgb(156pt)=(0.3176,0.7994,0.4975); rgb(157pt)=(0.3301,0.8002,0.4876); rgb(158pt)=(0.3424,0.8009,0.4774); rgb(159pt)=(0.3548,0.8016,0.4669); rgb(160pt)=(0.3671,0.8021,0.4563); rgb(161pt)=(0.3795,0.8026,0.4454); rgb(162pt)=(0.3921,0.8029,0.4344); rgb(163pt)=(0.405,0.8031,0.4233); rgb(164pt)=(0.4184,0.803,0.4122); rgb(165pt)=(0.4322,0.8028,0.4013); rgb(166pt)=(0.4463,0.8024,0.3904); rgb(167pt)=(0.4608,0.8018,0.3797); rgb(168pt)=(0.4753,0.8011,0.3691); rgb(169pt)=(0.4899,0.8002,0.3586); rgb(170pt)=(0.5044,0.7993,0.348); rgb(171pt)=(0.5187,0.7982,0.3374); rgb(172pt)=(0.5329,0.797,0.3267); rgb(173pt)=(0.547,0.7957,0.3159); rgb(175pt)=(0.5748,0.7929,0.2941); rgb(176pt)=(0.5886,0.7913,0.2833); rgb(177pt)=(0.6024,0.7896,0.2726); rgb(178pt)=(0.6161,0.7878,0.2622); rgb(179pt)=(0.6297,0.7859,0.2521); rgb(180pt)=(0.6433,0.7839,0.2423); rgb(181pt)=(0.6567,0.7818,0.2329); rgb(182pt)=(0.6701,0.7796,0.2239); rgb(183pt)=(0.6833,0.7773,0.2155); rgb(184pt)=(0.6963,0.775,0.2075); rgb(185pt)=(0.7091,0.7727,0.1998); rgb(186pt)=(0.7218,0.7703,0.1924); rgb(187pt)=(0.7344,0.7679,0.1852); rgb(188pt)=(0.7468,0.7654,0.1782); rgb(189pt)=(0.759,0.7629,0.1717); rgb(190pt)=(0.771,0.7604,0.1658); rgb(191pt)=(0.7829,0.7579,0.1608); rgb(192pt)=(0.7945,0.7554,0.157); rgb(193pt)=(0.806,0.7529,0.1546); rgb(194pt)=(0.8172,0.7505,0.1535); rgb(195pt)=(0.8281,0.7481,0.1536); rgb(196pt)=(0.8389,0.7457,0.1546); rgb(197pt)=(0.8495,0.7435,0.1564); rgb(198pt)=(0.86,0.7413,0.1587); rgb(199pt)=(0.8703,0.7392,0.1615); rgb(200pt)=(0.8804,0.7372,0.165); rgb(201pt)=(0.8903,0.7353,0.1695); rgb(202pt)=(0.9,0.7336,0.1749); rgb(203pt)=(0.9093,0.7321,0.1815); rgb(204pt)=(0.9184,0.7308,0.189); rgb(205pt)=(0.9272,0.7298,0.1973); rgb(206pt)=(0.9357,0.729,0.2061); rgb(207pt)=(0.944,0.7285,0.2151); rgb(208pt)=(0.9523,0.7284,0.2237); rgb(209pt)=(0.9606,0.7285,0.2312); rgb(210pt)=(0.9689,0.7292,0.2373); rgb(211pt)=(0.977,0.7304,0.2418); rgb(212pt)=(0.9842,0.733,0.2446); rgb(213pt)=(0.99,0.7365,0.2429); rgb(214pt)=(0.9946,0.7407,0.2394); rgb(215pt)=(0.9966,0.7458,0.2351); rgb(216pt)=(0.9971,0.7513,0.2309); rgb(217pt)=(0.9972,0.7569,0.2267); rgb(218pt)=(0.9971,0.7626,0.2224); rgb(219pt)=(0.9969,0.7683,0.2181); rgb(220pt)=(0.9966,0.774,0.2138); rgb(221pt)=(0.9962,0.7798,0.2095); rgb(222pt)=(0.9957,0.7856,0.2053); rgb(223pt)=(0.9949,0.7915,0.2012); rgb(224pt)=(0.9938,0.7974,0.1974); rgb(225pt)=(0.9923,0.8034,0.1939); rgb(226pt)=(0.9906,0.8095,0.1906); rgb(227pt)=(0.9885,0.8156,0.1875); rgb(228pt)=(0.9861,0.8218,0.1846); rgb(229pt)=(0.9835,0.828,0.1817); rgb(230pt)=(0.9807,0.8342,0.1787); rgb(231pt)=(0.9778,0.8404,0.1757); rgb(232pt)=(0.9748,0.8467,0.1726); rgb(233pt)=(0.972,0.8529,0.1695); rgb(234pt)=(0.9694,0.8591,0.1665); rgb(235pt)=(0.9671,0.8654,0.1636); rgb(236pt)=(0.9651,0.8716,0.1608); rgb(237pt)=(0.9634,0.8778,0.1582); rgb(238pt)=(0.9619,0.884,0.1557); rgb(239pt)=(0.9608,0.8902,0.1532); rgb(240pt)=(0.9601,0.8963,0.1507); rgb(241pt)=(0.9596,0.9023,0.148); rgb(242pt)=(0.9595,0.9084,0.145); rgb(243pt)=(0.9597,0.9143,0.1418); rgb(244pt)=(0.9601,0.9203,0.1382); rgb(245pt)=(0.9608,0.9262,0.1344); rgb(246pt)=(0.9618,0.932,0.1304); rgb(247pt)=(0.9629,0.9379,0.1261); rgb(248pt)=(0.9642,0.9437,0.1216); rgb(249pt)=(0.9657,0.9494,0.1168); rgb(250pt)=(0.9674,0.9552,0.1116); rgb(251pt)=(0.9692,0.9609,0.1061); rgb(252pt)=(0.9711,0.9667,0.1001); rgb(253pt)=(0.973,0.9724,0.0938); rgb(254pt)=(0.9749,0.9782,0.0872); rgb(255pt)=(0.9769,0.9839,0.0805)},
colorbar,
colorbar/width=0.1in,
colorbar style={ytick={0,0.25,0.5,0.75,1}, scaled y ticks = false, font=\tiny, yticklabel style={/pgf/number format/.cd,fixed,}, at={(5.25in,1in)},axis line style={very thin}} 
]
\addplot [forget plot] graphics [xmin=0.5, xmax=57.5, ymin=0.5, ymax=57.5] {figures/#2-2.png};
\end{axis}
\end{tikzpicture}
}

\newcommand{\myResPlotb}[2]{
\scriptsize
\begin{tikzpicture}
\begin{axis}[%
width=1in,
height=1in,
at={(0in,0in)},
scale only axis,
point meta min=0,
point meta max=#1,
axis on top,
separate axis lines,
every outer x axis line/.append style={black},
every x tick label/.append style={font=\color{black}},
every x tick/.append style={black},
xmin=0.5,
xmax=57.5,
xlabel = {Original},
xlabel near ticks,
every outer y axis line/.append style={black},
every y tick label/.append style={font=\color{black}},
every y tick/.append style={black},
y dir=reverse,
ymin=0.5,
ymax=57.5,
axis line style={draw=none},
ticks=none,
legend style={legend cell align=left, align=left, draw=black},
]
\addplot [forget plot] graphics [xmin=0.5, xmax=57.5, ymin=0.5, ymax=57.5] {figures/#2-1.png};
\end{axis}

\begin{axis}[%
width=1in,
height=1in,
at={(1.075in,0in)},
scale only axis,
point meta min=0,
point meta max=#1,
axis on top,
separate axis lines,
every outer x axis line/.append style={black},
every x tick label/.append style={font=\color{black}},
every x tick/.append style={black},
xmin=0.5,
xmax=57.5,
every outer y axis line/.append style={black},
every y tick label/.append style={font=\color{black}},
every y tick/.append style={black},
y dir=reverse,
ymin=0.5,
ymax=57.5,
xlabel = {SKA (ST)},
xlabel near ticks,
axis line style={draw=none},
ticks=none,
legend style={legend cell align=left, align=left, draw=black},
]
\addplot [forget plot] graphics [xmin=0.5, xmax=57.5, ymin=0.5, ymax=57.5] {figures/#2-4.png};
\end{axis}

\begin{axis}[%
width=1in,
height=1in,
at={(2.15in,0in)},
scale only axis,
point meta min=0,
point meta max=#1,
axis on top,
separate axis lines,
every outer x axis line/.append style={black},
every x tick label/.append style={font=\color{black}},
every x tick/.append style={black},
xmin=0.5,
xmax=57.5,
every outer y axis line/.append style={black},
every y tick label/.append style={font=\color{black}},
every y tick/.append style={black},
y dir=reverse,
ymin=0.5,
ymax=57.5,
xlabel = {FISTA (ST)},
xlabel near ticks,
axis line style={draw=none},
ticks=none,
legend style={legend cell align=left, align=left, draw=black},
]
\addplot [forget plot] graphics [xmin=0.5, xmax=57.5, ymin=0.5, ymax=57.5] {figures/#2-5.png};
\end{axis}

\begin{axis}[%
width=1in,
height=1in,
at={(3.225in,0in)},
scale only axis,
point meta min=0,
point meta max=#1,
axis on top,
separate axis lines,
every outer x axis line/.append style={black},
every x tick label/.append style={font=\color{black}},
every x tick/.append style={black},
xmin=0.5,
xmax=57.5,
every outer y axis line/.append style={black},
every y tick label/.append style={font=\color{black}},
every y tick/.append style={black},
y dir=reverse,
ymin=0.5,
ymax=57.5,
xlabel = {FISTA (NNG)},
xlabel near ticks,
axis line style={draw=none},
ticks=none,
legend style={legend cell align=left, align=left, draw=black},
]
\addplot [forget plot] graphics [xmin=0.5, xmax=57.5, ymin=0.5, ymax=57.5] {figures/#2-3.png};
\end{axis}

\begin{axis}[%
width=1in,
height=1in,
at={(4.3in,0in)},
scale only axis,
point meta min=0,
point meta max=#1,
axis on top,
separate axis lines,
every outer x axis line/.append style={black},
every x tick label/.append style={font=\color{black}},
every x tick/.append style={black},
xmin=0.5,
xmax=57.5,
xlabel = {SKA (NNG)},
xlabel near ticks,
every outer y axis line/.append style={black},
every y tick label/.append style={font=\color{black}},
every y tick/.append style={black},
y dir=reverse,
ymin=0.5,
ymax=57.5,
axis line style={draw=none},
ticks=none,
legend style={legend cell align=left, align=left, draw=black},
colormap={mymap}{[1pt] rgb(0pt)=(0.2422,0.1504,0.6603); rgb(1pt)=(0.2444,0.1534,0.6728); rgb(2pt)=(0.2464,0.1569,0.6847); rgb(3pt)=(0.2484,0.1607,0.6961); rgb(4pt)=(0.2503,0.1648,0.7071); rgb(5pt)=(0.2522,0.1689,0.7179); rgb(6pt)=(0.254,0.1732,0.7286); rgb(7pt)=(0.2558,0.1773,0.7393); rgb(8pt)=(0.2576,0.1814,0.7501); rgb(9pt)=(0.2594,0.1854,0.761); rgb(11pt)=(0.2628,0.1932,0.7828); rgb(12pt)=(0.2645,0.1972,0.7937); rgb(13pt)=(0.2661,0.2011,0.8043); rgb(14pt)=(0.2676,0.2052,0.8148); rgb(15pt)=(0.2691,0.2094,0.8249); rgb(16pt)=(0.2704,0.2138,0.8346); rgb(17pt)=(0.2717,0.2184,0.8439); rgb(18pt)=(0.2729,0.2231,0.8528); rgb(19pt)=(0.274,0.228,0.8612); rgb(20pt)=(0.2749,0.233,0.8692); rgb(21pt)=(0.2758,0.2382,0.8767); rgb(22pt)=(0.2766,0.2435,0.884); rgb(23pt)=(0.2774,0.2489,0.8908); rgb(24pt)=(0.2781,0.2543,0.8973); rgb(25pt)=(0.2788,0.2598,0.9035); rgb(26pt)=(0.2794,0.2653,0.9094); rgb(27pt)=(0.2798,0.2708,0.915); rgb(28pt)=(0.2802,0.2764,0.9204); rgb(29pt)=(0.2806,0.2819,0.9255); rgb(30pt)=(0.2809,0.2875,0.9305); rgb(31pt)=(0.2811,0.293,0.9352); rgb(32pt)=(0.2813,0.2985,0.9397); rgb(33pt)=(0.2814,0.304,0.9441); rgb(34pt)=(0.2814,0.3095,0.9483); rgb(35pt)=(0.2813,0.315,0.9524); rgb(36pt)=(0.2811,0.3204,0.9563); rgb(37pt)=(0.2809,0.3259,0.96); rgb(38pt)=(0.2807,0.3313,0.9636); rgb(39pt)=(0.2803,0.3367,0.967); rgb(40pt)=(0.2798,0.3421,0.9702); rgb(41pt)=(0.2791,0.3475,0.9733); rgb(42pt)=(0.2784,0.3529,0.9763); rgb(43pt)=(0.2776,0.3583,0.9791); rgb(44pt)=(0.2766,0.3638,0.9817); rgb(45pt)=(0.2754,0.3693,0.984); rgb(46pt)=(0.2741,0.3748,0.9862); rgb(47pt)=(0.2726,0.3804,0.9881); rgb(48pt)=(0.271,0.386,0.9898); rgb(49pt)=(0.2691,0.3916,0.9912); rgb(50pt)=(0.267,0.3973,0.9924); rgb(51pt)=(0.2647,0.403,0.9935); rgb(52pt)=(0.2621,0.4088,0.9946); rgb(53pt)=(0.2591,0.4145,0.9955); rgb(54pt)=(0.2556,0.4203,0.9965); rgb(55pt)=(0.2517,0.4261,0.9974); rgb(56pt)=(0.2473,0.4319,0.9983); rgb(57pt)=(0.2424,0.4378,0.9991); rgb(58pt)=(0.2369,0.4437,0.9996); rgb(59pt)=(0.2311,0.4497,0.9995); rgb(60pt)=(0.225,0.4559,0.9985); rgb(61pt)=(0.2189,0.462,0.9968); rgb(62pt)=(0.2128,0.4682,0.9948); rgb(63pt)=(0.2066,0.4743,0.9926); rgb(64pt)=(0.2006,0.4803,0.9906); rgb(65pt)=(0.195,0.4861,0.9887); rgb(66pt)=(0.1903,0.4919,0.9867); rgb(67pt)=(0.1869,0.4975,0.9844); rgb(68pt)=(0.1847,0.503,0.9819); rgb(69pt)=(0.1831,0.5084,0.9793); rgb(70pt)=(0.1818,0.5138,0.9766); rgb(71pt)=(0.1806,0.5191,0.9738); rgb(72pt)=(0.1795,0.5244,0.9709); rgb(73pt)=(0.1785,0.5296,0.9677); rgb(74pt)=(0.1778,0.5349,0.9641); rgb(75pt)=(0.1773,0.5401,0.9602); rgb(76pt)=(0.1768,0.5452,0.956); rgb(77pt)=(0.1764,0.5504,0.9516); rgb(78pt)=(0.1755,0.5554,0.9473); rgb(79pt)=(0.174,0.5605,0.9432); rgb(80pt)=(0.1716,0.5655,0.9393); rgb(81pt)=(0.1686,0.5705,0.9357); rgb(82pt)=(0.1649,0.5755,0.9323); rgb(83pt)=(0.161,0.5805,0.9289); rgb(84pt)=(0.1573,0.5854,0.9254); rgb(85pt)=(0.154,0.5902,0.9218); rgb(86pt)=(0.1513,0.595,0.9182); rgb(87pt)=(0.1492,0.5997,0.9147); rgb(88pt)=(0.1475,0.6043,0.9113); rgb(89pt)=(0.1461,0.6089,0.908); rgb(90pt)=(0.1446,0.6135,0.905); rgb(91pt)=(0.1429,0.618,0.9022); rgb(92pt)=(0.1408,0.6226,0.8998); rgb(93pt)=(0.1383,0.6272,0.8975); rgb(94pt)=(0.1354,0.6317,0.8953); rgb(95pt)=(0.1321,0.6363,0.8932); rgb(96pt)=(0.1288,0.6408,0.891); rgb(97pt)=(0.1253,0.6453,0.8887); rgb(98pt)=(0.1219,0.6497,0.8862); rgb(99pt)=(0.1185,0.6541,0.8834); rgb(100pt)=(0.1152,0.6584,0.8804); rgb(101pt)=(0.1119,0.6627,0.877); rgb(102pt)=(0.1085,0.6669,0.8734); rgb(103pt)=(0.1048,0.671,0.8695); rgb(104pt)=(0.1009,0.675,0.8653); rgb(105pt)=(0.0964,0.6789,0.8609); rgb(106pt)=(0.0914,0.6828,0.8562); rgb(107pt)=(0.0855,0.6865,0.8513); rgb(108pt)=(0.0789,0.6902,0.8462); rgb(109pt)=(0.0713,0.6938,0.8409); rgb(110pt)=(0.0628,0.6972,0.8355); rgb(111pt)=(0.0535,0.7006,0.8299); rgb(112pt)=(0.0433,0.7039,0.8242); rgb(113pt)=(0.0328,0.7071,0.8183); rgb(114pt)=(0.0234,0.7103,0.8124); rgb(115pt)=(0.0155,0.7133,0.8064); rgb(116pt)=(0.0091,0.7163,0.8003); rgb(117pt)=(0.0046,0.7192,0.7941); rgb(118pt)=(0.0019,0.722,0.7878); rgb(119pt)=(0.0009,0.7248,0.7815); rgb(120pt)=(0.0018,0.7275,0.7752); rgb(121pt)=(0.0046,0.7301,0.7688); rgb(122pt)=(0.0094,0.7327,0.7623); rgb(123pt)=(0.0162,0.7352,0.7558); rgb(124pt)=(0.0253,0.7376,0.7492); rgb(125pt)=(0.0369,0.74,0.7426); rgb(126pt)=(0.0504,0.7423,0.7359); rgb(127pt)=(0.0638,0.7446,0.7292); rgb(128pt)=(0.077,0.7468,0.7224); rgb(129pt)=(0.0899,0.7489,0.7156); rgb(130pt)=(0.1023,0.751,0.7088); rgb(131pt)=(0.1141,0.7531,0.7019); rgb(132pt)=(0.1252,0.7552,0.695); rgb(133pt)=(0.1354,0.7572,0.6881); rgb(134pt)=(0.1448,0.7593,0.6812); rgb(135pt)=(0.1532,0.7614,0.6741); rgb(136pt)=(0.1609,0.7635,0.6671); rgb(137pt)=(0.1678,0.7656,0.6599); rgb(138pt)=(0.1741,0.7678,0.6527); rgb(139pt)=(0.1799,0.7699,0.6454); rgb(140pt)=(0.1853,0.7721,0.6379); rgb(141pt)=(0.1905,0.7743,0.6303); rgb(142pt)=(0.1954,0.7765,0.6225); rgb(143pt)=(0.2003,0.7787,0.6146); rgb(144pt)=(0.2061,0.7808,0.6065); rgb(145pt)=(0.2118,0.7828,0.5983); rgb(146pt)=(0.2178,0.7849,0.5899); rgb(147pt)=(0.2244,0.7869,0.5813); rgb(148pt)=(0.2318,0.7887,0.5725); rgb(149pt)=(0.2401,0.7905,0.5636); rgb(150pt)=(0.2491,0.7922,0.5546); rgb(151pt)=(0.2589,0.7937,0.5454); rgb(152pt)=(0.2695,0.7951,0.536); rgb(153pt)=(0.2809,0.7964,0.5266); rgb(154pt)=(0.2929,0.7975,0.517); rgb(155pt)=(0.3052,0.7985,0.5074); rgb(156pt)=(0.3176,0.7994,0.4975); rgb(157pt)=(0.3301,0.8002,0.4876); rgb(158pt)=(0.3424,0.8009,0.4774); rgb(159pt)=(0.3548,0.8016,0.4669); rgb(160pt)=(0.3671,0.8021,0.4563); rgb(161pt)=(0.3795,0.8026,0.4454); rgb(162pt)=(0.3921,0.8029,0.4344); rgb(163pt)=(0.405,0.8031,0.4233); rgb(164pt)=(0.4184,0.803,0.4122); rgb(165pt)=(0.4322,0.8028,0.4013); rgb(166pt)=(0.4463,0.8024,0.3904); rgb(167pt)=(0.4608,0.8018,0.3797); rgb(168pt)=(0.4753,0.8011,0.3691); rgb(169pt)=(0.4899,0.8002,0.3586); rgb(170pt)=(0.5044,0.7993,0.348); rgb(171pt)=(0.5187,0.7982,0.3374); rgb(172pt)=(0.5329,0.797,0.3267); rgb(173pt)=(0.547,0.7957,0.3159); rgb(175pt)=(0.5748,0.7929,0.2941); rgb(176pt)=(0.5886,0.7913,0.2833); rgb(177pt)=(0.6024,0.7896,0.2726); rgb(178pt)=(0.6161,0.7878,0.2622); rgb(179pt)=(0.6297,0.7859,0.2521); rgb(180pt)=(0.6433,0.7839,0.2423); rgb(181pt)=(0.6567,0.7818,0.2329); rgb(182pt)=(0.6701,0.7796,0.2239); rgb(183pt)=(0.6833,0.7773,0.2155); rgb(184pt)=(0.6963,0.775,0.2075); rgb(185pt)=(0.7091,0.7727,0.1998); rgb(186pt)=(0.7218,0.7703,0.1924); rgb(187pt)=(0.7344,0.7679,0.1852); rgb(188pt)=(0.7468,0.7654,0.1782); rgb(189pt)=(0.759,0.7629,0.1717); rgb(190pt)=(0.771,0.7604,0.1658); rgb(191pt)=(0.7829,0.7579,0.1608); rgb(192pt)=(0.7945,0.7554,0.157); rgb(193pt)=(0.806,0.7529,0.1546); rgb(194pt)=(0.8172,0.7505,0.1535); rgb(195pt)=(0.8281,0.7481,0.1536); rgb(196pt)=(0.8389,0.7457,0.1546); rgb(197pt)=(0.8495,0.7435,0.1564); rgb(198pt)=(0.86,0.7413,0.1587); rgb(199pt)=(0.8703,0.7392,0.1615); rgb(200pt)=(0.8804,0.7372,0.165); rgb(201pt)=(0.8903,0.7353,0.1695); rgb(202pt)=(0.9,0.7336,0.1749); rgb(203pt)=(0.9093,0.7321,0.1815); rgb(204pt)=(0.9184,0.7308,0.189); rgb(205pt)=(0.9272,0.7298,0.1973); rgb(206pt)=(0.9357,0.729,0.2061); rgb(207pt)=(0.944,0.7285,0.2151); rgb(208pt)=(0.9523,0.7284,0.2237); rgb(209pt)=(0.9606,0.7285,0.2312); rgb(210pt)=(0.9689,0.7292,0.2373); rgb(211pt)=(0.977,0.7304,0.2418); rgb(212pt)=(0.9842,0.733,0.2446); rgb(213pt)=(0.99,0.7365,0.2429); rgb(214pt)=(0.9946,0.7407,0.2394); rgb(215pt)=(0.9966,0.7458,0.2351); rgb(216pt)=(0.9971,0.7513,0.2309); rgb(217pt)=(0.9972,0.7569,0.2267); rgb(218pt)=(0.9971,0.7626,0.2224); rgb(219pt)=(0.9969,0.7683,0.2181); rgb(220pt)=(0.9966,0.774,0.2138); rgb(221pt)=(0.9962,0.7798,0.2095); rgb(222pt)=(0.9957,0.7856,0.2053); rgb(223pt)=(0.9949,0.7915,0.2012); rgb(224pt)=(0.9938,0.7974,0.1974); rgb(225pt)=(0.9923,0.8034,0.1939); rgb(226pt)=(0.9906,0.8095,0.1906); rgb(227pt)=(0.9885,0.8156,0.1875); rgb(228pt)=(0.9861,0.8218,0.1846); rgb(229pt)=(0.9835,0.828,0.1817); rgb(230pt)=(0.9807,0.8342,0.1787); rgb(231pt)=(0.9778,0.8404,0.1757); rgb(232pt)=(0.9748,0.8467,0.1726); rgb(233pt)=(0.972,0.8529,0.1695); rgb(234pt)=(0.9694,0.8591,0.1665); rgb(235pt)=(0.9671,0.8654,0.1636); rgb(236pt)=(0.9651,0.8716,0.1608); rgb(237pt)=(0.9634,0.8778,0.1582); rgb(238pt)=(0.9619,0.884,0.1557); rgb(239pt)=(0.9608,0.8902,0.1532); rgb(240pt)=(0.9601,0.8963,0.1507); rgb(241pt)=(0.9596,0.9023,0.148); rgb(242pt)=(0.9595,0.9084,0.145); rgb(243pt)=(0.9597,0.9143,0.1418); rgb(244pt)=(0.9601,0.9203,0.1382); rgb(245pt)=(0.9608,0.9262,0.1344); rgb(246pt)=(0.9618,0.932,0.1304); rgb(247pt)=(0.9629,0.9379,0.1261); rgb(248pt)=(0.9642,0.9437,0.1216); rgb(249pt)=(0.9657,0.9494,0.1168); rgb(250pt)=(0.9674,0.9552,0.1116); rgb(251pt)=(0.9692,0.9609,0.1061); rgb(252pt)=(0.9711,0.9667,0.1001); rgb(253pt)=(0.973,0.9724,0.0938); rgb(254pt)=(0.9749,0.9782,0.0872); rgb(255pt)=(0.9769,0.9839,0.0805)},
colorbar,
colorbar/width=0.1in,
colorbar style={ytick={0,0.25,0.5,0.75,1.00}, scaled y ticks = false, font=\tiny, yticklabel style={/pgf/number format/.cd,fixed,}, at={(5.25in,1in)},axis line style={very thin}} 
]
\addplot [forget plot] graphics [xmin=0.5, xmax=57.5, ymin=0.5, ymax=57.5] {figures/#2-2.png};
\end{axis}
\end{tikzpicture}
}

\begin{figure*}
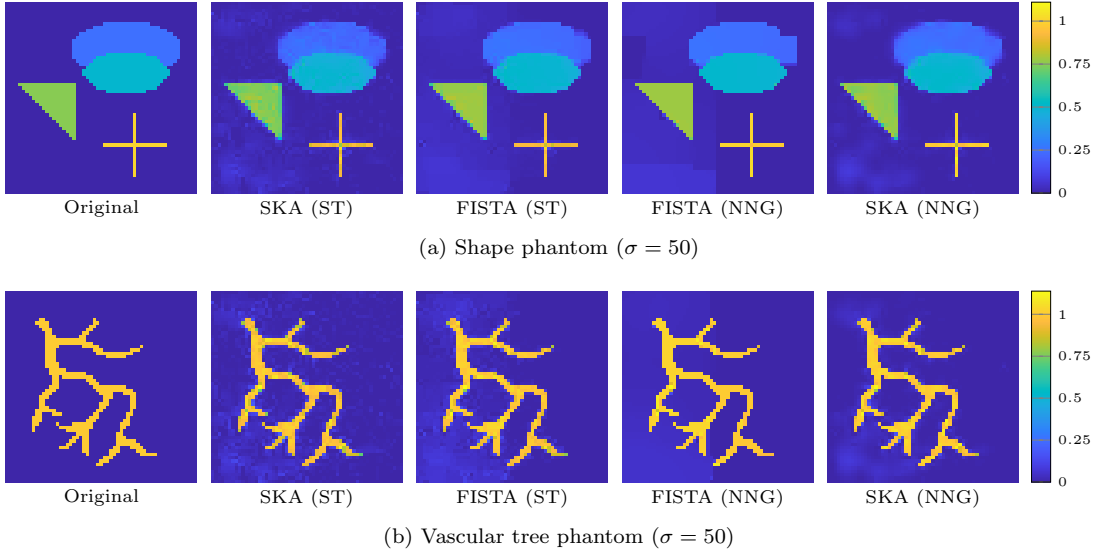
%
\centering
	\subfloat[Shape phantom ($\sigma=50$)]{
		\hspace{-6.5pt}\myResPlot{1.11}{resBlurC}
	}\\
	\subfloat[Vascular tree phantom ($\sigma=50$)]{
		\hspace{-6.5pt}\myResPlotb{1.138}{resVascC}
	}
\caption{Reconstructed tracer concentration for both phantoms at a large noise level. Shown are the original phantom (left) and the reconstructed images for the four algorithmic approaches. For each phantom the corresponding colormap is normalized and valid for all five images.}%
\label{fig:resC}%
\end{figure*}

\subsection{Simulated Data}

The PSNR, SSIM and computation times (averaged over 50 independent runs) for the two phantoms are summarized in the first four rows in Table \ref{tab:1} for a low noise level of $\sigma=10$ and a large noise level $\sigma=50$ for the FISTA and the proposed SKA method based on the soft-thresholding (ST) and the non-negative Garrote (NNG). The SKA method with NNG thresholding gives the best PSNR results for all test cases. In addition, it has significantly better SSIM values and reduces the computation times by at least a factor of 2.5 compared to the classical FISTA. The speed of convergence of the proposed Kaczmarz based algorithm in comparison with FISTA is visualized in Fig.~\ref{fig:res1} for all four test scenarios. With the chosen stopping criterion in \eqref{eq:conv}, convergence is reached for all algorithms, while our Kaczmarz based approaches reaches convergence faster than its FISTA counterpart. In more detail, the Kaczmarz based approach with soft-thresholding leads to a poorer image quality while showing slightly faster computation times to reach the stopping criterion.

The reconstructed tracer concentration in the case of large noise ($\sigma=50$) is visualized in Fig.~\ref{fig:resC}. It shows that the algorithms based on NNG thresholding leads to a smoother particle concentration at less background noise compared to soft-thresholding, where the Kaczmarz based approach yields the superior reconstruction performance. This can be quantified by comparing the actual reconstructed values for certain areas in the shape phantom. The triangle part has an original value of $0.75$, reconstruction using FISTA (NNG) gives $0.7761\pm0.0083$ and with soft-thresholding $0.7665\pm0.0281$, while our SKA (NNG) approach yields $0.7592\pm0.0292$ and with soft-thresholding $0.733\pm0.0382$. The background of the shape phantom (original value of 0) is reconstructed to the following values: FISTA (NNG) $0.0186\pm0.0228$, FISTA (ST) $0.0211\pm0.0224$, SKA (NNG) $0.0094\pm0.0175$ and SKA (ST) $0.0134\pm0.0193$.

In summary, the proposed SKA with NNG thresholding shows the best performance in terms of image quality. In addition, it also shows faster convergence, leading to consistently faster computation times. In the following section, we are going to put these results into perspective with state-of-the-art regularization algorithms currently used in the context of MPI.

\subsection{Comparison with state-of-the-art Algorithms}

There are currently two state-of-the-art algorithms for image reconstruction in MPI. The first one, is the Tikhonov regularized Kaczmarz (reg.~KA) algorithm proposed in \cite{Knopp2010}, which solves the regularization problem in \eqref{eq:minprob}. The second one solves the optimization problem \eqref{eq:fl}. Corresponding Fused Lasso algorithms are introduced in \cite{Storath2017} and \cite{ilbey2017} which both solve the same minimization problem, and, have quite similar performance in terms of PSNR as well as computation times \cite[Fig.~3-6]{ilbey2017}. 

For each of the two algorithms, the corresponding regularization parameters are chosen such that the PSNR is maximized (with the same stopping criterion as given in \eqref{eq:conv}). The results are summarized in the bottom two rows of Table~\ref{tab:1}. As expected, the regularized Kaczmarz approach gives the worst PSNR. Storath et al.'s Fused Lasso leads to much better reconstruction results, improving PSNR as well as SSIM at the cost of increasing computation time. It is outperformed, however, by the proposed SKA (NNG) approach, in terms of PSNR and SSIM and most significantly in terms of computation times.   
This implies that our proposed optimization problem in \eqref{eq:minpr} is better suited to reconstruct the particle concentration in MPI than the problems \eqref{eq:minprob} or \eqref{eq:fl}.

\subsection{Experimental Data}

\begin{figure}[t!]%
\centering
\subfloat[Cone phantom (taken from \cite{Knopp2019})]{
	\includegraphics[width=0.7\columnwidth]{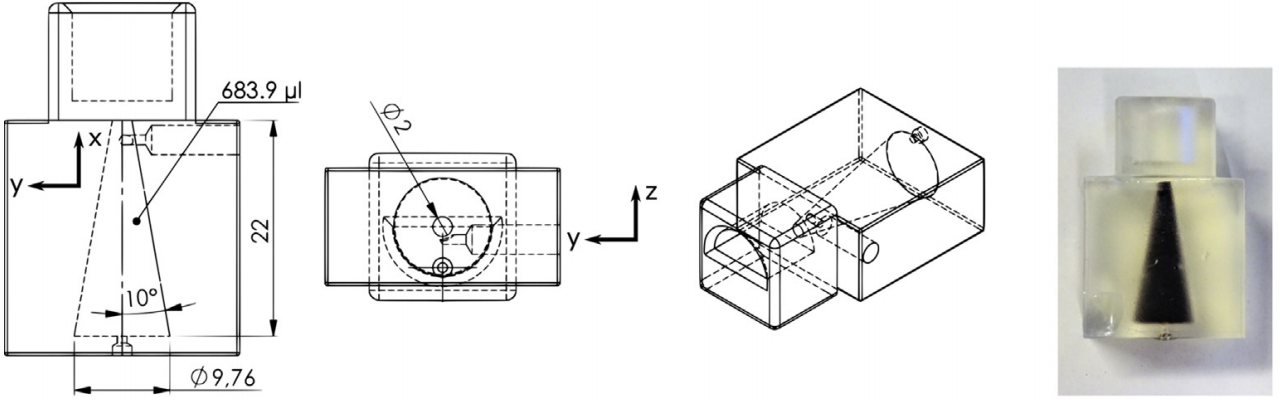}
	\label{fig:shaph}
}\\
\subfloat[Resolution phantom (taken from \cite{Knopp2019})]{
	\includegraphics[width=0.7\columnwidth]{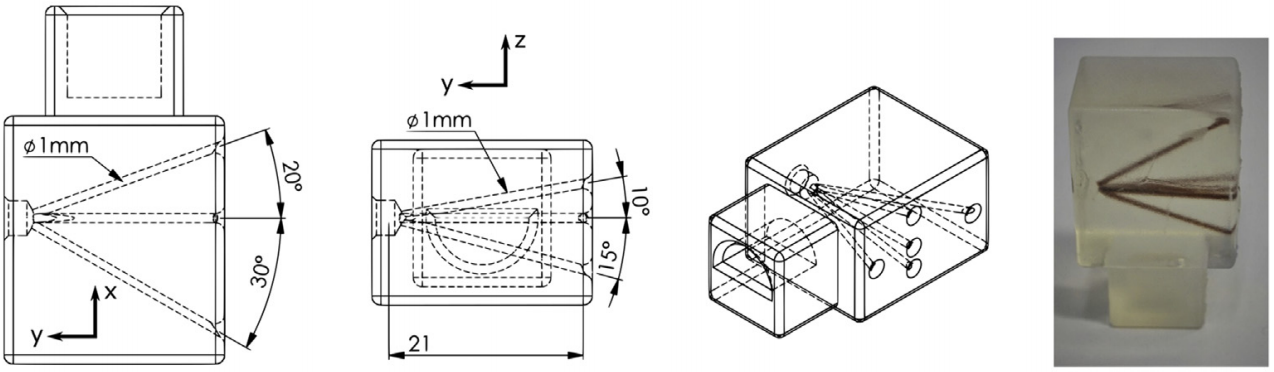}
	\label{fig:resph}
}
\caption{Schematics of the used phantoms. All images are taken from \cite{Knopp2019}.}%
\label{fig:phantoms}%
\end{figure}

The openly available 3D data set used in this work is acquired using a Bruker preclinical MPI scanner \cite{Knopp2019}. The 3D system matrix is calibrated with the tracer perimag (100\,mmol/l) and a delta sample of size 2\,mm\,$\times$\,2\,mm\,$\times$1\,mm. The drive-field amplitude is 12\,mT\,$\times$\,12\,mT\,$\times$\,12\,mT at a frequency of 2.5\,MHz/102\,$\times$\,2.5\,Mhz/96\,$\times$\,2.5\,Mhz/99 and a selection-field gradient of -1.0\,T/m\,$\times$\,-1.0\,T/m\,$\times$\,2.0\,T/m. The FOV is discretized on a grid of size 37\,$\times$\,37\,$\times$\,37 amounting to 50,653 voxel. The signal for each voxel is averaged over 1,000 measurements and Fourier transformed resulting in a complex-valued system matrix $A$ of size 80,787\,$\times$\,50,653. During the system matrix acquisition, background measurements are taken after every 37th voxel scan. The system matrix is then background corrected by interpolating the resulting background measurements and subtracting them from the system matrix \cite{knopp2019correction}. 
The phantoms used are the "cone" and the "resolution" phantom which are shown in Fig.~\ref{fig:resph} and \ref{fig:shaph} taken from \cite{Knopp2019}.
Both phantoms are filled with perimag at a concentration of 50\,mmol/l. The measurement vector $b$ then results by taking the mean over 1,000 phantom measurements and subtracting the mean over the same amount of empty scanner measurements.

The parameter settings are as follows. For reconstruction we consider only frequencies which are larger than 70\,kHz and have an SNR larger than 3. This reduces the number of matrix rows from 80,787 down to 3,271. Hence, the resulting system matrix $A$ has size 3,271\,$\times$\,50,653. For FISTA the largest eigenvalue of $A^*A$ is needed, which implies an eigenvalue decomposition of a matrix of size 50,653\,$\times$\,50,653. In our case we used the iterative Arnoldi's method to approximate the largest eigenvalue. With ten iterations of Arnoldi's method the computation of the largest eigenvalue takes about 6 minutes and is not included in the subsequent evaluation of computation times. Note, that whenever the parameter setting is changed, i.e., SNR or frequency threshold, this also affects the eigenvalue of $A^*A$.

The reconstruction of the tracer concentration is based on a three dimensional undecimated wavelet transform as described in Section \ref{sec:swt}. For the Kaczmarz-based approaches, the stopping criterion is set to $\epsilon_r = 10^{-2}$, while both FISTA approaches need $\epsilon_r = 10^{-3}$ in order to reach a reconstruction quality which is comparable with the Kaczmarz-based approaches. The chosen regularization parameters $\lambda$  are summarized in Table \ref{tab:2} for each of the two phantoms. For comparison, we include the results for the regularized Kaczmarz as the state-of-the-art reconstruction method. Storath et al.'s Fused Lasso approach was not able to handle the large size of the system matrix and is therefore omitted in the following. Note that computation times of Storath et al.'s Fused Lasso would be significant larger, mainly due to the discretization scheme of the total variation norm.

\begin{table}%
\centering
\caption{Regularization parameters and computation times for the reconstruction results of the cone and resolution phantom shown in Fig.~\ref{fig:resOP}. Computation times are averaged over 20 runs.}
\label{tab:2}
\vspace{10pt}
\begin{tabular}{@{}rrrrr@{}} \toprule
& \multicolumn{2}{c}{cone ph.} & \multicolumn{2}{c}{resolution ph.} \\
 \cmidrule(r){2-3} \cmidrule(r){4-5} 
& $\lambda$ & time & $\lambda$ & time \\\midrule
SKA (NNG)   & 2.5e-2 & 14.9  & 7e-4 & 33.5  \\
SKA (ST)    & 1e-2   & 11.4  & 3e-4 & 26.3  \\
FISTA (NNG) & 3e-4   & 158.7 & 5e-5 & 305.6 \\
FISTA (ST)  & 2e-4   & 142.6 & 1e-5 & 240.8 \\
Reg.~KA		  & 1e-2   & 11.6  & 1e-4 & 22.8  \\
\bottomrule
\end{tabular}
\end{table}

\newcommand{\myThreeDplotb}[3]{
\begin{tikzpicture}
\scriptsize
\begin{axis}[%
width=2.0cm,
height=2.0cm,
at={(0cm,2.8cm)},
scale only axis,
point meta min=0,
point meta max=#1,
axis on top,
separate axis lines,
every outer x axis line/.append style={black},
every x tick label/.append style={font=\color{black}},
every x tick/.append style={black},
xmin=0.5,
xmax=19.5,
xtick={\empty},
every outer y axis line/.append style={black},
every y tick label/.append style={font=\color{black}},
every y tick/.append style={black},
y dir=reverse,
ymin=0.5,
ymax=19.5,
ytick={\empty},
axis line style={draw=none},
ticks=none,
title style={font=\scriptsize,yshift=-7pt},
title={z=19},
legend style={legend cell align=left, align=left, draw=black}
]
\addplot [forget plot] graphics [xmin=0.5, xmax=19.5, ymin=0.5, ymax=19.5] {figures/#2-1.png};
\end{axis}

\begin{axis}[%
width=2.0cm,
height=1.0cm,
at={(0cm,1.40cm)},
scale only axis,
point meta min=0,
point meta max=#1,
axis on top,
separate axis lines,
xmin=0.5,
xmax=19.5,
xtick={\empty},
ymin=0.5,
ymax=19.5,
ytick={\empty},
axis line style={draw=none},
ticks=none,
title style={font=\scriptsize,yshift=-8pt},
title={y=19},
legend style={legend cell align=left, align=left, draw=black}
]
\addplot [forget plot] graphics [xmin=0.5, xmax=19.5, ymin=0.5, ymax=19.5] {figures/#2-2.png};
\end{axis}

\begin{axis}[%
width=2.0cm,
height=1.0cm,
at={(0cm,0cm)},
scale only axis,
point meta min=0,
point meta max=#1,
axis on top,
separate axis lines,
every outer x axis line/.append style={black},
every x tick label/.append style={font=\color{black}},
every x tick/.append style={black},
xmin=0.5,
xmax=19.5,
xtick={\empty},
xlabel={#3},
xlabel style={yshift=10pt},
every outer y axis line/.append style={black},
every y tick label/.append style={font=\color{black}},
every y tick/.append style={black},
y dir=reverse,
ymin=0.5,
ymax=19.5,
ytick={\empty},
axis line style={draw=none},
ticks=none,
title style={font=\scriptsize,yshift=-7pt},
title={x=19},
legend style={legend cell align=left, align=left, draw=black},
colormap={mymap}{[1pt] rgb(0pt)=(0.2422,0.1504,0.6603); rgb(1pt)=(0.25039,0.164995,0.707614); rgb(2pt)=(0.257771,0.181781,0.751138); rgb(3pt)=(0.264729,0.197757,0.795214); rgb(4pt)=(0.270648,0.214676,0.836371); rgb(5pt)=(0.275114,0.234238,0.870986); rgb(6pt)=(0.2783,0.255871,0.899071); rgb(7pt)=(0.280333,0.278233,0.9221); rgb(8pt)=(0.281338,0.300595,0.941376); rgb(9pt)=(0.281014,0.322757,0.957886); rgb(10pt)=(0.279467,0.344671,0.971676); rgb(11pt)=(0.275971,0.366681,0.982905); rgb(12pt)=(0.269914,0.3892,0.9906); rgb(13pt)=(0.260243,0.412329,0.995157); rgb(14pt)=(0.244033,0.435833,0.998833); rgb(15pt)=(0.220643,0.460257,0.997286); rgb(16pt)=(0.196333,0.484719,0.989152); rgb(17pt)=(0.183405,0.507371,0.979795); rgb(18pt)=(0.178643,0.528857,0.968157); rgb(19pt)=(0.176438,0.549905,0.952019); rgb(20pt)=(0.168743,0.570262,0.935871); rgb(21pt)=(0.154,0.5902,0.9218); rgb(22pt)=(0.146029,0.609119,0.907857); rgb(23pt)=(0.138024,0.627629,0.89729); rgb(24pt)=(0.124814,0.645929,0.888343); rgb(25pt)=(0.111252,0.6635,0.876314); rgb(26pt)=(0.0952095,0.679829,0.859781); rgb(27pt)=(0.0688714,0.694771,0.839357); rgb(28pt)=(0.0296667,0.708167,0.816333); rgb(29pt)=(0.00357143,0.720267,0.7917); rgb(30pt)=(0.00665714,0.731214,0.766014); rgb(31pt)=(0.0433286,0.741095,0.73941); rgb(32pt)=(0.0963952,0.75,0.712038); rgb(33pt)=(0.140771,0.7584,0.684157); rgb(34pt)=(0.1717,0.766962,0.655443); rgb(35pt)=(0.193767,0.775767,0.6251); rgb(36pt)=(0.216086,0.7843,0.5923); rgb(37pt)=(0.246957,0.791795,0.556743); rgb(38pt)=(0.290614,0.79729,0.518829); rgb(39pt)=(0.340643,0.8008,0.478857); rgb(40pt)=(0.3909,0.802871,0.435448); rgb(41pt)=(0.445629,0.802419,0.390919); rgb(42pt)=(0.5044,0.7993,0.348); rgb(43pt)=(0.561562,0.794233,0.304481); rgb(44pt)=(0.617395,0.787619,0.261238); rgb(45pt)=(0.671986,0.779271,0.2227); rgb(46pt)=(0.7242,0.769843,0.191029); rgb(47pt)=(0.773833,0.759805,0.16461); rgb(48pt)=(0.820314,0.749814,0.153529); rgb(49pt)=(0.863433,0.7406,0.159633); rgb(50pt)=(0.903543,0.733029,0.177414); rgb(51pt)=(0.939257,0.728786,0.209957); rgb(52pt)=(0.972757,0.729771,0.239443); rgb(53pt)=(0.995648,0.743371,0.237148); rgb(54pt)=(0.996986,0.765857,0.219943); rgb(55pt)=(0.995205,0.789252,0.202762); rgb(56pt)=(0.9892,0.813567,0.188533); rgb(57pt)=(0.978629,0.838629,0.176557); rgb(58pt)=(0.967648,0.8639,0.16429); rgb(59pt)=(0.96101,0.889019,0.153676); rgb(60pt)=(0.959671,0.913457,0.142257); rgb(61pt)=(0.962795,0.937338,0.12651); rgb(62pt)=(0.969114,0.960629,0.106362); rgb(63pt)=(0.9769,0.9839,0.0805)},
colorbar,
colorbar style={at={(1.85cm,4.8cm)},width=0.25cm,height=4.8cm,axis line style={very thin}, } 
]
\addplot [forget plot] graphics [xmin=0.5, xmax=19.5, ymin=0.5, ymax=19.5] {figures/#2-3.png};
\end{axis}
\end{tikzpicture}%
}

\newcommand{\myThreeDplota}[3]{
\begin{tikzpicture}
\scriptsize
\begin{axis}[%
width=2.5cm,
height=2.5cm,
at={(0cm,3.3cm)},
scale only axis,
point meta min=0,
point meta max=#1,
axis on top,
separate axis lines,
xmin=0.5,
xmax=19.5,
xtick={\empty},
y dir=reverse,
ymin=0.5,
ymax=19.5,
ytick={\empty},
axis line style={draw=none},
ticks=none,
title style={font=\scriptsize,yshift=-7pt},
title={z=19},
legend style={legend cell align=left, align=left, draw=black}
]
\addplot [forget plot] graphics [xmin=0.5, xmax=19.5, ymin=0.5, ymax=19.5] {figures/#2-1.png};
\end{axis}

\begin{axis}[%
width=2.5cm,
height=1.25cm,
at={(0cm,1.65cm)},
scale only axis,
point meta min=0,
point meta max=#1,
axis on top,
separate axis lines,
xmin=0.5,
xmax=19.5,
xtick={\empty},
ymin=0.5,
ymax=19.5,
ytick={\empty},
axis line style={draw=none},
ticks=none,
title style={font=\scriptsize,yshift=-8pt},
title={y=19},
legend style={legend cell align=left, align=left, draw=black}
]
\addplot [forget plot] graphics [xmin=0.5, xmax=19.5, ymin=0.5, ymax=19.5] {figures/#2-2.png};
\end{axis}

\begin{axis}[%
width=2.5cm,
height=1.25cm,
at={(0cm,0cm)},
scale only axis,
point meta min=0,
point meta max=#1,
axis on top,
separate axis lines,
xmin=0.5,
xmax=19.5,
xtick={\empty},
xlabel={#3},
xlabel style={yshift=10pt},
y dir=reverse,
ymin=0.5,
ymax=19.5,
ytick={\empty},
axis line style={draw=none},
ticks=none,
title style={font=\scriptsize,yshift=-7pt},
title={x=19},
]
\addplot [forget plot] graphics [xmin=0.5, xmax=19.5, ymin=0.5, ymax=19.5] {figures/#2-3.png};
\end{axis}
\end{tikzpicture}%
}

\begin{figure*}[t!]
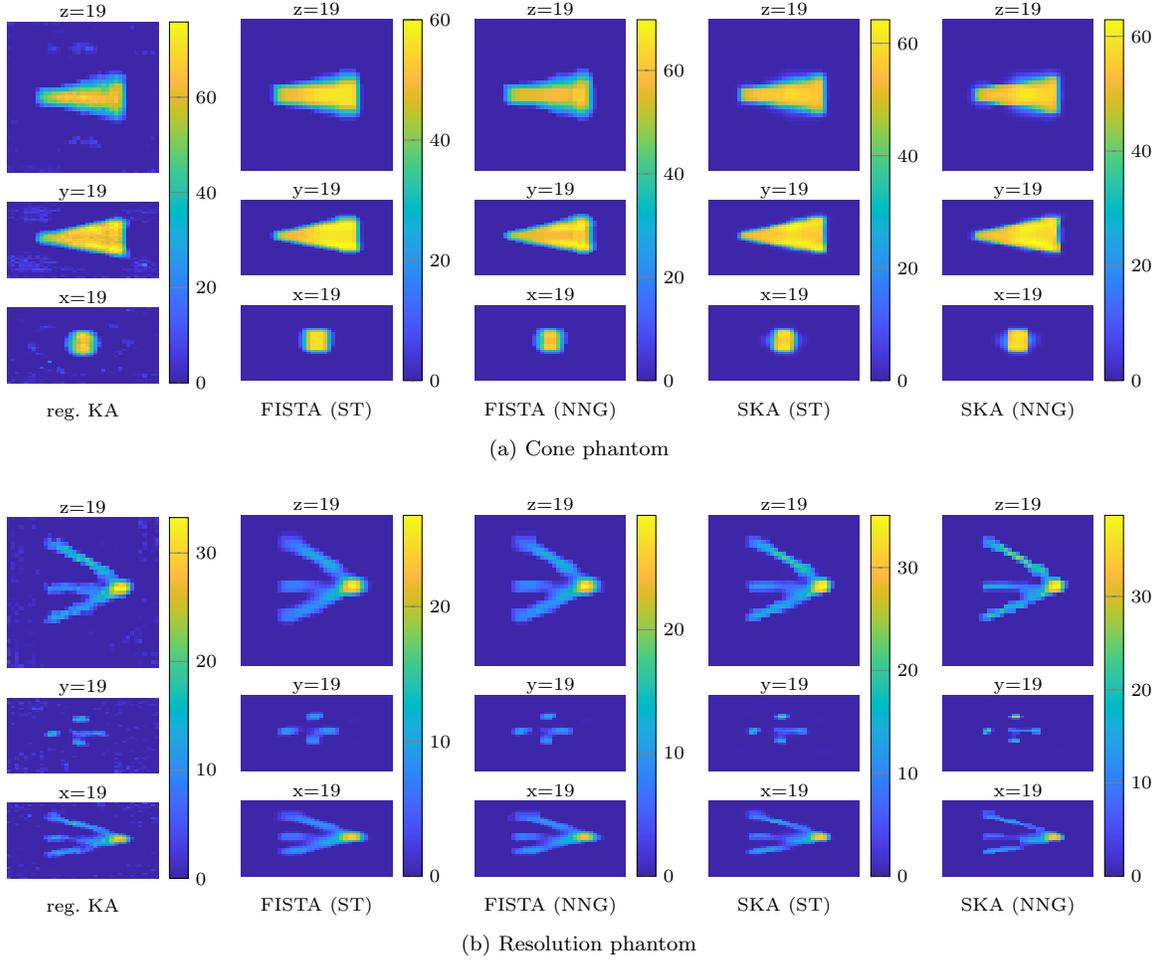
%
\centering
\subfloat[Cone phantom]{
\myThreeDplotb{75.80}{ShapeP_HighRes_regKA}{reg.~KA}
\myThreeDplotb{60.01}{ResP_HighRes_fistast}{FISTA (ST)}
\myThreeDplotb{69.81}{ResP_HighRes_fistaew}{FISTA (NNG)}
\myThreeDplotb{64.18}{ShapeP_HighRes_skast}{SKA (ST)}
\myThreeDplotb{62.9}{ShapeP_HighRes_skaew}{SKA (NNG)}
\label{fig:resOpA}
}\\
\subfloat[Resolution phantom]{
\myThreeDplotb{33.26}{ResolutionP_HighRes_regKA}{reg.~KA}
\myThreeDplotb{26.77}{ResolutionP_HighRes_fistast}{FISTA (ST)}
\myThreeDplotb{29.29}{ResolutionP_HighRes_fistaew}{FISTA (NNG)}
\myThreeDplotb{35.08}{ResolutionP_HighRes_skast1}{SKA (ST)}
\myThreeDplotb{38.70}{ResolutionP_HighRes_skaew1}{SKA (NNG)}
\label{fig:resOpB}
}
\caption{Reconstructed 3D tracer concentration in layer view for the cone and resolution phantom. }%
\label{fig:resOP}%
\end{figure*}

The resulting tracer reconstructions are visualized in Figure \ref{fig:resOP}. The three dimensional cone in Fig.~\ref{fig:resOpA} is shown in layer view, i.e., three 2D images, showing a triangle in the $xy$\nobreakdash- and $xz$\nobreakdash-plane and a circle in the $yz$\nobreakdash-plane. Whereas the regularized Kaczmarz still has artefacts in the background, there are none in all wavelet based approaches. The computation times for the cone reconstruction are summarized in Table \ref{tab:2} and clearly favor the proposed Kaczmarz based approaches, accelerating reconstruction by a factor of 14 compared to FISTA. 
The reconstruction results of the resolution phantom are visualized in Fig.~\ref{fig:resOpB}. Again, the wavelet based approaches show no background artefacts compared to the regularized Kaczmarz. Compared to FISTA, the sparse Kaczmarz approaches lead to better reconstructions in terms of edge reconstruction, e.g. compare the layer-view in the $xy$\nobreakdash-plane of the SKA (NNG) and FISTA (NNG) reconstructions. This is mainly due to insufficient iteration steps for the FISTA based approaches. The reconstruction results would improve further by increasing the number of iterations, i.e., decreasing the stopping criterion. In contrast, the sparse Kaczmarz approaches only need about 26 and 34 seconds with soft-thresholding and NNG-thresholding respectively. Interestingly, the fastest reconstruction results from the reg.~KA approach. This can be explained by quite similar iteration numbers compared to SKA (ST) to reach the stopping criterion. Note that this is quite sensitive to the chosen SNR threshold used in the preprocessing of the system matrix.


\section{Discussion and Conclusion}
In this work we introduced a novel row-action algorithm based on sparsity in the wavelet domain for image reconstruction in magnetic particle imaging. Our proposed sparse Kaczmarz method is compared to a classical FISTA approach, where both algorithms solve the same minimization problem. The performance evaluation on simulated MPI data shows a significant improvement compared to FISTA and current state-of-the-art algorithms in terms of PSNR, SSIM but also in terms of computation time. All algorithms were also evaluated on a large openly available MPI data set. The proposed sparse Kaczmarz approaches yield a reasonable reconstruction quality at a fraction of the computation time of the FISTA approaches. It thus combines the advantages of the row-action methods with the more sophisticated image priors of Storath et al.'s Fused Lasso approach: a superior reconstructed image quality with less noise, better edge preservation and a considerably faster convergence rate. In addition, the sparse Kaczmarz only requires a single row of the system matrix at each inner iteration. This reduces the computational overhead of computing $A^*A$ and its largest eigenvalue and makes it possible to handle system matrices that do not fit into computer memory.  For potential usages in real-time MPI applications such as angioplasty and real-time catheter tracking but also for offline MPI reconstruction, the proposed approach would be a notable benefit.

In contrast to the state-of-the-art Fused Lasso algorithm, the proposed method only has a single regularization parameter $\lambda$ instead of two in the Fused Lasso case, which simplifies fine-tuning. The tuning of $\lambda$, however, directly affects the reconstructed image quality and is therefore a crucial task. An in-depth analysis of different parameter choices such as the discrepancy principle, quasi-optimality criterion \cite{Kluth2019}, the Akaike Information Criteria (AIC) \cite[Ch.~7]{Hastie2009} or a simple wavelet based noise estimator \cite{Mallat2009} is outside the scope of this manuscript and leaves room for further improvement. For large system matrices it might be useful to apply a block Kaczmarz approach \cite{Needell2015} or a weighted block approach \cite{Necoara2019}, which can further accelerate convergence. A different approach to further speed up the computation time is an asynchronous implementation of the Kaczmarz iteration \cite{Liu2014}. But also different thresholding operators, such as the Persistent Empirical Wiener (PEW) in \cite{Kowalski2014} performing the thresholding not based on a single pixel but on a neighborhood of surrounding pixels, as well as the whitening approaches described in \cite{Kluth2019} can further improve the reconstructed tracer concentration. 



\section*{Acknowledgment}
The authors would like to thank A.~von Gladi{\ss} for providing the hybrid system matrix, used in the simulation study.
This work was supported by the German Federal Ministry of Education and Research (BMBF grant number 05M16WFA).

\bibliographystyle{IEEEtran}

\begin{thebibliography}{10}
\providecommand{\url}[1]{#1}
\csname url@samestyle\endcsname
\providecommand{\newblock}{\relax}
\providecommand{\bibinfo}[2]{#2}
\providecommand{\BIBentrySTDinterwordspacing}{\spaceskip=0pt\relax}
\providecommand{\BIBentryALTinterwordstretchfactor}{4}
\providecommand{\BIBentryALTinterwordspacing}{\spaceskip=\fontdimen2\font plus
\BIBentryALTinterwordstretchfactor\fontdimen3\font minus
  \fontdimen4\font\relax}
\providecommand{\BIBforeignlanguage}[2]{{%
\expandafter\ifx\csname l@#1\endcsname\relax
\typeout{** WARNING: IEEEtran.bst: No hyphenation pattern has been}%
\typeout{** loaded for the language `#1'. Using the pattern for}%
\typeout{** the default language instead.}%
\else
\language=\csname l@#1\endcsname
\fi
#2}}
\providecommand{\BIBdecl}{\relax}
\BIBdecl

\bibitem{Gleich2005}
B.~Gleich and J.~Weizenecker, ``Tomographic imaging using the nonlinear
  response of magnetic particles,'' \emph{Nature}, vol. 435, no. 7046, pp.
  1214--1217, jun 2005.

\bibitem{Panagiotopoulos2015}
N.~Panagiotopoulos, F.~Vogt, J.~Barkhausen, T.~M. Buzug, R.~L. Duschka,
  K.~Lüdtke-Buzug, M.~Ahlborg, G.~Bringout, C.~Debbeler, M.~Gräser,
  C.~Kaethner, J.~Stelzner, H.~Medimagh, and J.~Haegele, ``Magnetic particle
  imaging: current developments and future directions,'' \emph{International
  Journal of Nanomedicine}, p. 3097, apr 2015.

\bibitem{Bakenecker2018}
A.~Bakenecker, M.~Ahlborg, C.~Debbeler, C.~Kaethner, and K.~Lüdtke-Buzug,
  ``Magnetic particle imaging,'' in \emph{Precision Medicine}.\hskip 1em plus
  0.5em minus 0.4em\relax Elsevier, 2018, pp. 183--228.

\bibitem{Vaalma2017}
S.~Vaalma, J.~Rahmer, N.~Panagiotopoulos, R.~L. Duschka, J.~Borgert,
  J.~Barkhausen, F.~M. Vogt, and J.~Haegele, ``Magnetic particle imaging
  ({MPI}): Experimental quantification of vascular stenosis using stationary
  stenosis phantoms,'' \emph{{PLOS} {ONE}}, vol.~12, no.~1, p. e0168902, jan
  2017.

\bibitem{Yu2017}
E.~Y. Yu, M.~Bishop, B.~Zheng, R.~M. Ferguson, A.~P. Khandhar, S.~J. Kemp,
  K.~M. Krishnan, P.~W. Goodwill, and S.~M. Conolly, ``Magnetic particle
  imaging: A novel in vivo imaging platform for cancer detection,'' \emph{Nano
  Letters}, vol.~17, no.~3, pp. 1648--1654, feb 2017.

\bibitem{Tay2018}
Z.~W. Tay, P.~Chandrasekharan, A.~Chiu-Lam, D.~W. Hensley, R.~Dhavalikar, X.~Y.
  Zhou, E.~Y. Yu, P.~W. Goodwill, B.~Zheng, C.~Rinaldi, and S.~M. Conolly,
  ``Magnetic particle imaging-guided heating in vivo using gradient fields for
  arbitrary localization of magnetic hyperthermia therapy,'' \emph{{ACS} Nano},
  vol.~12, no.~4, pp. 3699--3713, mar 2018.

\bibitem{Le2017}
T.-A. Le, X.~Zhang, A.~K. Hoshiar, and J.~Yoon, ``Real-time two-dimensional
  magnetic particle imaging for electromagnetic navigation in targeted drug
  delivery,'' \emph{Sensors}, vol.~17, no.~9, p. 2050, sep 2017.

\bibitem{Knopp2008}
T.~Knopp, S.~Biederer, T.~Sattel, J.~Weizenecker, B.~Gleich, J.~Borgert, and
  T.~M. Buzug, ``Trajectory analysis for magnetic particle imaging,''
  \emph{Physics in Medicine and Biology}, vol.~54, no.~2, pp. 385--397, dec
  2008.

\bibitem{Graeser2017}
M.~Graeser, A.~von Gladiss, M.~Weber, and T.~M. Buzug, ``Two dimensional
  magnetic particle spectrometry,'' \emph{Physics in Medicine and Biology},
  vol.~62, no.~9, pp. 3378--3391, apr 2017.

\bibitem{ozaslan2019fully}
A.~Ozaslan, A.~Alacaoglu, O.~Demirel, T.~{\c{C}}ukur, and E.~Saritas, ``Fully
  automated gridding reconstruction for non-cartesian x-space magnetic particle
  imaging,'' \emph{Physics in Medicine \& Biology}, vol.~64, no.~16, p. 165018,
  2019.

\bibitem{Knopp2010}
T.~Knopp, J.~Rahmer, T.~F. Sattel, S.~Biederer, J.~Weizenecker, B.~Gleich,
  J.~Borgert, and T.~M. Buzug, ``Weighted iterative reconstruction for magnetic
  particle imaging,'' \emph{Physics in Medicine and Biology}, vol.~55, no.~6,
  pp. 1577--1589, feb 2010.

\bibitem{Kluth2019}
T.~Kluth and B.~Jin, ``Enhanced reconstruction in magnetic particle imaging by
  whitening and randomized {SVD} approximation,'' \emph{Physics in Medicine
  {\&} Biology}, vol.~64, no.~12, p. 125026, jun 2019.

\bibitem{Knopp2016a}
T.~Knopp and M.~Hofmann, ``Online reconstruction of 3d magnetic particle
  imaging data,'' \emph{Physics in Medicine and Biology}, vol.~61, no.~11, pp.
  N257--N267, may 2016.

\bibitem{Moeddel2019}
M.~Möddel, T.~Knopp, R.~Werner, D.~Weller, and J.~M. Salamon, ``Toward
  employing the full potential of magnetic particle imaging: exploring
  visualization techniques and clinical use cases for real-time 3d vascular
  imaging,'' in \emph{Medical Imaging 2019: Biomedical Applications in
  Molecular, Structural, and Functional Imaging}, B.~Gimi and A.~Krol,
  Eds.\hskip 1em plus 0.5em minus 0.4em\relax {SPIE}, mar 2019.

\bibitem{Storath2017}
M.~Storath, C.~Brandt, M.~Hofmann, T.~Knopp, J.~Salamon, A.~Weber, and
  A.~Weinmann, ``Edge preserving and noise reducing reconstruction for magnetic
  particle imaging,'' \emph{{IEEE} Transactions on Medical Imaging}, vol.~36,
  no.~1, pp. 74--85, jan 2017.

\bibitem{ilbey2017}
S.~Ilbey, C.~Top, A.~Güngör, T.~Cukur, E.~Saritas, and H.~Guven, ``Comparison
  of system-matrix-based and projection-based reconstructions for field free
  line magnetic particle imaging,'' \emph{International Journal on Magnetic
  Particle Imaging}, vol.~3, 2017.

\bibitem{Knopp2019}
\BIBentryALTinterwordspacing
T.~Knopp, P.~Szwargulski, F.~Griese, and M.~Gräser, ``Openmpidata: An
  initiative for freely accessible magnetic particle imaging data,'' \emph{Data
  in Brief}, p. 104971, 2019. [Online]. Available:
  \url{http://www.sciencedirect.com/science/article/pii/S2352340919313265}
\BIBentrySTDinterwordspacing

\bibitem{Fadili2009}
M.~Fadili and J.-L. Starck, ``Monotone operator splitting for optimization
  problems in sparse recovery,'' in \emph{2009 16th {IEEE} International
  Conference on Image Processing ({ICIP})}.\hskip 1em plus 0.5em minus
  0.4em\relax {IEEE}, nov 2009.

\bibitem{Kowalski2014}
M.~Kowalski, ``Thresholding {RULES} and iterative shrinkage/thresholding
  algorithm: A convergence study,'' in \emph{2014 {IEEE} International
  Conference on Image Processing ({ICIP})}.\hskip 1em plus 0.5em minus
  0.4em\relax {IEEE}, oct 2014.

\bibitem{lustig2007sparse}
M.~Lustig, D.~Donoho, and J.~M. Pauly, ``Sparse mri: The application of
  compressed sensing for rapid mr imaging,'' \emph{Magnetic Resonance in
  Medicine: An Official Journal of the International Society for Magnetic
  Resonance in Medicine}, vol.~58, no.~6, pp. 1182--1195, 2007.

\bibitem{mallat89}
S.~Mallat, ``A theory for multiresolution signal decomposition: the wavelet
  representation,'' \emph{{IEEE} Transactions on Pattern Analysis and Machine
  Intelligence}, vol.~11, no.~7, pp. 674--693, jul 1989.

\bibitem{Nason1995}
G.~P. Nason and B.~W. Silverman, ``The stationary wavelet transform and some
  statistical applications,'' in \emph{Wavelets and Statistics}.\hskip 1em plus
  0.5em minus 0.4em\relax Springer New York, 1995, pp. 281--299.

\bibitem{Kaczmarz37}
S.~Kaczmarz, ``Angen{\"a}herte {Aufl{\"o}sung} von {Systemen} linearer
  {Gleichungen},'' \emph{Bull. Internat. Acad. Polon. Sci. Lettres A}, 1937.

\bibitem{Schoepfer2018}
F.~Schöpfer and D.~A. Lorenz, ``Linear convergence of the randomized sparse
  kaczmarz method,'' \emph{Mathematical Programming}, vol. 173, no. 1-2, pp.
  509--536, jan 2018.

\bibitem{Strohmer2008}
T.~Strohmer and R.~Vershynin, ``A randomized kaczmarz algorithm with
  exponential convergence,'' \emph{Journal of Fourier Analysis and
  Applications}, vol.~15, no.~2, pp. 262--278, apr 2008.

\bibitem{Lorenz2014a}
D.~A. Lorenz, F.~Schöpfer, and S.~Wenger, ``The linearized bregman method via
  split feasibility problems: Analysis and generalizations,'' \emph{{SIAM}
  Journal on Imaging Sciences}, vol.~7, no.~2, pp. 1237--1262, jan 2014.

\bibitem{Trussell1985}
H.~Trussell and M.~Civanlar, ``The landweber iteration and projection onto
  convex sets,'' \emph{{IEEE} Transactions on Acoustics, Speech, and Signal
  Processing}, vol.~33, no.~6, pp. 1632--1634, dec 1985.

\bibitem{Beck2009a}
A.~Beck and M.~Teboulle, ``A fast iterative shrinkage-thresholding algorithm
  for linear inverse problems,'' \emph{{SIAM} Journal on Imaging Sciences},
  vol.~2, no.~1, pp. 183--202, jan 2009.

\bibitem{Boyd2010}
S.~Boyd, ``Distributed optimization and statistical learning via the
  alternating direction method of multipliers,'' \emph{Foundations and
  Trends{\textregistered} in Machine Learning}, vol.~3, no.~1, pp. 1--122,
  2010.

\bibitem{Combettes2011}
P.~L. Combettes and J.-C. Pesquet, ``Proximal splitting methods in signal
  processing,'' in \emph{Springer Optimization and Its Applications}.\hskip 1em
  plus 0.5em minus 0.4em\relax Springer New York, 2011, pp. 185--212.

\bibitem{Chen2018}
X.~Chen, ``The kaczmarz algorithm, row action methods, and statistical learning
  algorithms,'' pp. 115--127, 2018.

\bibitem{Scherzer2007}
O.~Scherzer, A.~Leit{\~{a}}o, and M.~Haltmeier, ``Kaczmarz methods for
  regularizing nonlinear ill-posed equations i: convergence analysis,''
  \emph{Inverse Problems and Imaging}, vol.~1, no.~2, pp. 289--298, apr 2007.

\bibitem{Bayram2016}
I.~Bayram, ``On the convergence of the iterative shrinkage/thresholding
  algorithm with a weakly convex penalty,'' \emph{{IEEE} Transactions on Signal
  Processing}, vol.~64, no.~6, pp. 1597--1608, mar 2016.

\bibitem{Guo2017}
K.~Guo, D.~Han, and X.~Yuan, ``Convergence analysis of douglas--rachford
  splitting method for {\textquotedblleft}strongly + weakly{\textquotedblright}
  convex programming,'' \emph{{SIAM} Journal on Numerical Analysis}, vol.~55,
  no.~4, pp. 1549--1577, jan 2017.

\bibitem{Zhang2019}
T.~Zhang and Z.~Shen, ``A fundamental proof of convergence of alternating
  direction method of multipliers for weakly convex optimization,''
  \emph{Journal of Inequalities and Applications}, vol. 2019, no.~1, may 2019.

\bibitem{Gladiss2017}
A.~von Gladiss, M.~Graeser, P.~Szwargulski, T.~Knopp, and T.~M. Buzug, ``Hybrid
  system calibration for multidimensional magnetic particle imaging,''
  \emph{Physics in Medicine and Biology}, vol.~62, no.~9, pp. 3392--3406, apr
  2017.

\bibitem{Bathke2017}
C.~Bathke, T.~Kluth, C.~Brandt, and P.~Maass,
  ``\BIBforeignlanguage{en}{Improved image reconstruction in magnetic particle
  imaging using structural a priori information},''
  \emph{\BIBforeignlanguage{en}{International Journal on Magnetic Particle
  Imaging}}, vol. Vol.3, 2017.

\bibitem{Hansen1994}
P.~C. Hansen, ``{REGULARIZATION} {TOOLS}: A matlab package for analysis and
  solution of discrete ill-posed problems,'' \emph{Numerical Algorithms},
  vol.~6, no.~1, pp. 1--35, mar 1994.

\bibitem{Wang2004}
Z.~Wang, A.~Bovik, H.~Sheikh, and E.~Simoncelli, ``Image quality assessment:
  From error visibility to structural similarity,'' \emph{{IEEE} Transactions
  on Image Processing}, vol.~13, no.~4, pp. 600--612, apr 2004.

\bibitem{knopp2019correction}
T.~Knopp, N.~Gdaniec, R.~Rehr, M.~Gr{\"a}ser, and T.~Gerkmann, ``Correction of
  linear system drifts in magnetic particle imaging,'' \emph{Physics in
  Medicine \& Biology}, vol.~64, no.~12, p. 125013, 2019.

\bibitem{Hastie2009}
T.~Hastie, R.~Tibshirani, and J.~Friedman, \emph{The Elements of Statistical
  Learning}.\hskip 1em plus 0.5em minus 0.4em\relax Springer New York, 2009.

\bibitem{Mallat2009}
\BIBentryALTinterwordspacing
S.~Mallat, \emph{A Wavelet Tour of Signal Processing}.\hskip 1em plus 0.5em
  minus 0.4em\relax Elsevier LTD, Oxford, 2009. [Online]. Available:
  \url{https://www.ebook.de/de/product/7596640/stephane_mallat_a_wavelet_tour_of_signal_processing.html}
\BIBentrySTDinterwordspacing

\bibitem{Needell2015}
D.~Needell, R.~Zhao, and A.~Zouzias, ``Randomized block kaczmarz method with
  projection for solving least squares,'' \emph{Linear Algebra and its
  Applications}, vol. 484, pp. 322--343, nov 2015.

\bibitem{Necoara2019}
I.~Necoara, ``Faster randomized block kaczmarz algorithms,''
  \emph{arXiv:1902.09946}, 2019.

\bibitem{Liu2014}
J.~Liu, S.~J. Wright, and S.~Sridhar, ``An asynchronous parallel randomized
  kaczmarz algorithm,'' \emph{arXiv}, 2014.

\end{thebibliography}

\end{document}